\documentclass[11pt]{amsart}
\usepackage{xcolor}
\usepackage[utf8]{inputenc}
\usepackage{graphicx}
\usepackage{amssymb,amsmath,amsgen,amsopn,amsxtra,amsbsy,amscd}
\usepackage{xcolor}

\def\({\left(}
\def\){\right)}

\def\eb{\varepsilon}

\def\al{\alpha}
\def\rw{\rightarrow}
\def\Dx{\Delta}
\def\nx{\nabla}

\def\Om{\Omega}

\def\la{\lambda}

\def\R {\mathbb{R}}

\def\E {{\mathcal E}}

\def\LL {{\mathcal L}}

\def\M {{\mathcal M}}

\def\t{\tilde}

\def \p {\partial}

\def \and{\qquad\text{and}\qquad}

\def\Nx{\nabla}
\def \p {\partial}
\def\Dx{\Delta}
\newcommand{\be}{\begin{equation} }
\newcommand{\ee}{\end{equation} }

\newtheorem{proposition}{Proposition}[section]
\newtheorem{theorem}[proposition]{Theorem}

\theoremstyle{definition}

\newtheorem{remark}[proposition]{Remark}

\numberwithin{equation}{section}

\def \no#1#2#3 {{\bf #1} (#3), #2.}
\def \eds#1#2#3 {#1, #2, #3.}

\begin{document}
\title[Chevron pattern equations]{Chevron pattern equations: exponential attractor and global stabilization}

\author{H.~Kalantarova$^\ast$, V.~Kalantarov$^\dagger$ and O.~Vantzos$^\ddagger$}
\address{$^\ast$Department of Materials Science and Engineering, Technion-Israel Institute of Technology, 3200003 Haifa, Israel}
\address{$^\dagger$Department of Mathematics, Ko{\c{c}} University, Istanbul, Turkey}
\address{$^\dagger$Azerbaijan State Oil and Industry University, Baku, Azerbaijan }
\address{$^\ddagger$Lightricks Ltd., Jerusalem, Israel}


\date{\today}


\maketitle

\begin{abstract}
The initial boundary value problem for a nonlinear system of equations modeling the chevron patterns is studied in one and two spatial dimensions. The existence of an exponential attractor and the stabilization of the zero steady state solution through application of a finite-dimensional feedback control is proved in two spatial dimensions. The stabilization of an arbitrary fixed solution is shown in one spatial dimension along with relevant numerical results.
\end{abstract}

\section{Introduction}\label{s0}
We consider the following coupled system of equations introduced to model chevron patterns observed in the dielectric regime of electroconvection in nematic liquid crystals
\be\label{a0}
\tau \partial_{t}A=A+\Dx A - \phi^2A-|A|^2A-2i c_1\phi \partial_{y}A+i\beta A\partial_{y}\phi,
\ee
\be\label{a1}
\partial_{t}\phi=D_1\partial^{2}_{x}\phi+D_2\partial^{2}_{y}\phi-h\phi +\phi\lvert A\rvert^{2}-c_{2}\mbox{Im}\left[ A^{\ast}\partial_{y}A\right],
\ee
where  $\tau>0$, $D_1>0$, $D_2>0$, $c_1\geq0$, $c_2\geq0$, $h\geq0$, $\beta \in \R$ are given parameters, the complex valued function $A$ ($A^{\ast}$ denotes its complex conjugate), which describes the amplitude of the convection pattern and the real valued function $\phi$, which denotes the angle of the director from the x-axis are unknown functions. This system was proposed by Rossberg et al. \cite{Rossberg_diss}, \cite{Rossberg1996weakly} to model the formation and evolution of patterns. In his work, Rossberg shows that chevron patterns are observed in simulations of \eqref{a0}, \eqref{a1}. The application of liquid crystals in display devices, namely LCD screens and their use in electronics has made the study of the electro-optical effects in liquid crystals interesting both to researchers in academia and in industry.

The rest of the paper is organized as follows. Section \ref{sec:pr} is devoted to preliminary theorems and inequalities. In section \ref{sec:ea}, we prove the existence of an exponential attractor in two space dimensions. This is an improvement of a result from \cite{KKV}, where the existence of a global attractor for the system \eqref{a0}-\eqref{a1} is shown. In section \ref{sec:fs}, we use the finite-dimensional feedback control approach  for stabilizing the zero solution again in two space dimensions. Finally, in section 5 we study chevron patterns in one spatial dimension and prove global exponential stabilization to any fixed, not necessarily stationary, solution.

\section*{Notation}
Throughout this paper, $\mathbb{R}$ denotes the set of real numbers, $\Omega\subset\mathbb{R}^{2}$ denotes a bounded domain with sufficiently smooth boundary denoted by $\partial \Omega$.

\noindent $(\cdot,\cdot)$ and $\lVert\cdot\rVert$ denote the inner product and the norm induced by it in $L^{2}(\Omega)$, respectively. That is, for $f,g \in L^{2}(\Omega)$
\begin{equation}
(f,g):=\int_{\Omega}f(x,y)g^{\ast}(x,y)dxdy,\quad \lVert f\rVert:=\left(\int_{\Omega}\lvert f(x,y)\rvert^{2}dxdy\right)^{1/2}.\nonumber
\end{equation}
We are also using the following notation
$$
V^0:=L^2\times L^2,\quad V^1:=H^1_0\times H^1_0,\quad  V^2:=H^{2}\cap H^1_0\times H^2\cap H^1_0.
$$

\section{Preliminaries}\label{sec:pr}
In order to keep this work self contained, this section provides a theorem and inequalities that we need to prove results in the future sections.
 \begin{theorem}\label{exun} (\cite{KKV})
If $c_1\in [0,1)$ or $c_1\ge 2c_2$, $h>0$, the system of equations \eqref{a0}-\eqref{a1} together with the following initial and boundary conditions
\be\label{a2}
A\Big|_{t=0}=A_0, \ \ \phi\Big|_{t=0}=\phi_0, \ \ A\Big|_{\p \Om}=0, \ \ \phi\Big|_{\p \Om }=0,
\ee
where $A_{0},\phi_{0}\in L^{2}(\Omega)$, has a  unique  weak solution
\begin{equation}
A,\phi \in C([0,T]; L^2(\Om)) \cap L^2([0,T]; H_0^1(\Om)),\quad \forall T>0,  \ \
\end{equation}
such that
\be\label{aest1}
\|\phi(t)\|\le M_0, \ \ \|A(t)\|\le M_0, \ \ \forall t>0,
\ee
 and
\be\label{aest2}
\int_0^T\|\nx \phi(t)\|^{2}dt \le M_T, \ \  \int_0^T\|\nx A(t)\|^{2}dt\le M_T, \ \forall T>0,
\ee
where $M_{0}$ denotes a generic constant depending only on $\lVert A_{0}\rVert$, $\lVert\phi_{0}\rVert$ and $\lvert\Omega\rvert$, and where $M_{T}$ is a generic constant, which depends also on T, besides $\lVert A_{0}\rVert$, $\lVert\phi_{0}\rVert$ and $\lvert\Omega\rvert$. In other words this problem generates a compact semigroup $S(t$), $t\ge0,$ in the phase space $V^0$. Moreover  this semigroup has a global attractor that belongs to $V^1$.
\end{theorem}

\begin{itemize}
\item {\it Young\rq{}s inequality}
For each $a,b>0$ and  $\epsilon>0$
\be\label{Young}
ab \le \epsilon \frac{a^p}p+ \frac1{\epsilon^{q/p}}\frac{b^q}q,
\ee
where $p,q>0$ and $\frac1p+\frac1q=1.$

\item{\it}
For $u\in H^{2}(\Omega)\cap H^{1}_{0}(\Omega)$ the following estimates are valid
\begin{equation}\label{in:ell}
\nu_{0}\lVert\Delta u\rVert^{2}\leq \lVert\mathcal{L}u\rVert^{2}\leq\nu_{1}\lVert\Delta u\rVert^{2}
\end{equation}
where $\mathcal{L}$ is a second order uniformly elliptic operator.

\item {\it Poincar\'e-Friedrichs (P-F) inequality}
\be\label{PF}
\|u\|^2\le \la_1^{-1}\| \nabla u\|^2, \ \ \forall u\in  H_0^1(\Om)
\ee
and the inequality
\be\label{PFN}
\sum\limits_{k=N+1}^\infty|(u, w_k)|^2\le \la_{N+1}^{-1}\| \nabla u\|^2, \ \ \forall u\in  H_0^1(\Om),
\ee
where $0<\lambda_{1}\leq\lambda_{2}\leq\ldots\leq\lambda_{N+1}\leq\ldots$ are the eigenvalues of the Laplace operator, $-\Delta $,  under the homogeneous Dirichlet's boundary condition and $w_{k}$, for $k=1,2,\ldots$ are the corresponding eigenfunctions.
\item{\textit{Ladyzhenskaya inequality}}
\begin{equation}\label{ineq:lad}
\lVert  u\rVert_{L^{4}(\Omega)}\leq 2^{1/4}\lVert u\rVert^{1/2}\lVert \nabla u\rVert^{1/2}
\end{equation}
which is valid for $u\in H^{1}_{0}(\Omega)$ with $\Omega\subset\mathbb{R}^{2}$.
\item {\it The $1D$ Agmon  inequality}
\be\label{Ag}
\max_{x\in [0,L]}|u(x)|^2\le C_0 \|u\|\|\p_xu\|, \ \ \forall u\in H_0^1(0,L).
\ee
\item {\it The $2D$ Agmon  inequality}
\be\label{Ag2}
\max_{x\in \Om}|u(x)|^2\le C_0 \|u\|\|\Delta u\|, \ \ \forall u\in H^2(\Om) \cap H_0^1(\Om),
\ee
where $\Om \subset \R^2.$

\item {\it Gagliardo-Nirenberg inequality}

\noindent Let $\Omega\subset\mathbb{R}^{n}$ and $u\in W^{m,r}(\Omega)\cap L^{q}(\Omega)$, $1\leq r,q\leq\infty$. For any integer $j$, $0\leq j\leq m$ and for any number $\hat{a}$ in the interval $\frac{j}{m}\leq\hat{a}\leq 1$, set
\begin{equation}
\frac{1}{p}-\frac{j}{n}=\hat{a}\left(\frac{1}{r}-\frac{m}{n}\right)+(1-\hat{a})\frac{1}{q}.\nonumber
\end{equation}
If $m-j-\frac{n}{r}$ is a nonnegative integer, then
\begin{equation}\label{GN2}
\lVert D^{j}u\rVert_{L^{p}(\Omega)}\leq C\lVert u\rVert^{\hat{a}}_{W^{m,r}(\Omega)}\lVert u\rVert^{1-\hat{a}}_{L^{q}(\Omega)},
\end{equation}
where the constant $C$ depends only on $\Omega$, $r$, $m$, $q$, $j$ and $\hat{a}$.
\end{itemize}


\section{Exponential Attractor}\label{sec:ea}

In this section, we are going to show that the semigroup $S(t):V^0\rw V^0$, $t\ge 0$, associated to the problem \eqref{a0}-\eqref{a1} and \eqref{a2} has an exponential attractor, i.e., there exists a set
$\M \subset V^0$ which satisfies the following conditions:
\begin{enumerate}
\item $\M \subset V^0$ is compact and has a finite fractal dimension,
\item it is positively invariant, $S(t)\M\subset \M, \ \forall t>0$,
\item $\M$ attracts each bounded set $B\subset V^0$  with an exponential rate, i.e.,
for each bounded set $B\subset V^0$
$$
dist \ (S(t)B,\M)\le Q(\|B\|_{V^0})e^{-\alpha t}, \ t>0,
$$
\end{enumerate}
where the distance function is defined by
\begin{equation}
dist \ (S(t)B,\M):=\sup_{a\in S(t)B}\inf_{b\in\M}\lVert a-b\rVert_{V^{0}},\nonumber
\end{equation}
$\al>0$ and a monotone function $Q(\cdot)$ are independent of $B$.

\noindent For the literature about exponential attractors for various dissipative dynamical systems we refer to \cite{EFNT}, \cite{MEZ}, \cite{MiZe} and references therein. To prove the existence of an exponential attractor of the semigroup generated by the problem \eqref{a0}-\eqref{a2}, we use the following theorem
\begin{theorem}\label{MZ}(\cite{MiZe})
 Let $E$ and $E_1$ be two Banach spaces such that $E_{1}$ is compactly embedded in $E$. Assume that $S(t):E\to E$ is a semigroup, which possesses a compact  absorbing set $\mathcal B$ in $E$:
$$
S(t_*)\mathcal B\subset \mathcal B\ \mbox{ for }\ t_*>0.
$$
Assume further that,\\
a)	There exists $K>0$ such that
\begin{equation}\label{sma1}
\|S(t_{*})\xi_1-S(t_{*})\xi_2\|_{E_{1}}\le K\|\xi_1-\xi_2\|_E,\  \forall\xi_1,\xi_2\in\mathcal B.
\end{equation}
b)	The map $(t,\xi)\to S(t)\xi$ is H\"{o}lder continuous (or Lipschitz continuous) on $[0,t_*]\times\mathcal B$.

Then the semigroup $S(t)$ possesses an exponential attractor $\mathcal M$ in $E$, which is a subset of $\mathcal B$.
\end{theorem}

We proceed by showing that the semigroup $S(t):V^{0}\rightarrow V^{0}$ satisfies the conditions stated in Theorem \ref{MZ}. First, we show the existence of a compact absorbing set in $V^{0}$ for $S(t)$. We recall the following dissipative estimates
\begin{multline}\label{chD0}
\frac d{dt}\left[\tau\| A (t)\|^2 +\delta_{0}\|\phi(t)\|^2\right]
+k_0\left[\tau\| A (t)\|^2 +\delta_{0}\|\phi(t)\|^2\right]\\+ \delta_0
\left[\| \Nx A (t)\|^2 +D_{0}\|\Nx \phi(t)\|^2\right]\le |\Om|
\end{multline}
and
\begin{multline}\label{A12}
\frac{d}{dt}\left[\tau\|\nabla A(t)\|^{2}+(\LL \phi(t),\phi(t))\right]+\nu_{0}\|\Delta\phi\rVert^{2}+2^{-1} \|\Delta A\|^{2}\\
\leq C\left[\tau\|\nabla A\|^{2}\!+(\LL\phi,\phi)+1\right]\left[\tau\|\nabla A\|^{2}+ (\LL \phi,\phi)\right]
\end{multline}
which are derived in \cite[(10)]{KKV} and \cite[(38)]{KKV}, respectively. Here 
\begin{eqnarray}
\label{par1}&&k_0:=\min\{\tau^{-1},h\},\quad \delta_{0}:=2(1-c_{1})/(2+c_{2})\mbox{ where }c_{1}<1,\\
\label{par2}&&D_0:=\min\{D_1,D_2\},\quad -\mathcal{L}\phi:=-D_{1}\partial_{x}^{2}\phi-D_{2}\partial_{y}^{2}\phi,
\end{eqnarray}
and $C$ is a constant that depends on $\tau$, $\nu_{0}$, $c_{1}$, $c_{2}$, $\beta$, $D_{1}$, $D_{2}$ and $M_{0}$, which is defined in the statement of Theorem \ref{exun}. It follows from the estimate \eqref{chD0} that, the set
$$B_0:=\left\{[A,\phi]\subset V^0: \tau \|A\|^2+\delta_0\|\phi\|^2\le 2k_0^{-1}|\Om|\right\}$$
is a positively invariant absorbing set for the semigroup $S(t)$, $t\ge0$. The estimate \eqref{chD0} also implies that
\be\label{tt1}
\delta_{0}\int_t^{t+1}\left[\tau\|\nabla A(s)\|^{2}+(\LL \phi(s),\phi(s))\right]ds\le |\Om|(1+2k_0^{-1}),\quad \forall t>0.
\ee
Employing the uniform Gronwall lemma and the estimate \eqref{A12}, we deduce from \eqref{tt1} that
\be\label{chD2}
\tau\|\nabla A(t)\|^{2}+(\LL \phi(t),\phi(t))\le R_0,  \ \  \forall t\ge1,
\ee
where 
\begin{equation}
R_{0}=\frac{(1+C)\lvert\Omega\rvert(1+2k_{0}^{-1})}{\delta_{0}}e^{C\lvert\Omega\rvert(1+2k_{0}^{-1})/\delta_{0}}.\nonumber
\end{equation}
Hence $B_0\subset V^1$, and since $V^{1}$ is compactly embedded into $V^{0}$, $B_{0}$ is compact in $V^0$.

Next, we show that the condition a) is satisfied. Let $[A_0,\phi_0]$ and $[\t A_0,\t \phi_0]$ be arbitrary two elements of  $B_0$. Then $[a,\psi]=:[\tilde{A}-A,\tilde{\phi}-\phi]$ the  difference of $[A(t),\phi(t)]=S(t)[A_0,\phi_0]$ and $[\t A(t),\t \phi]= S(t)[\tilde{A}_{0},\tilde{\phi}_{0}]$ is a solution of the system
\begin{multline}\label{sm1}
\tau \p_t a =a+\Dx a-\phi^2a-(\tilde{\phi}^2-\phi^2)\t A -\lvert\tilde{A}\rvert^{2}\tilde{A}+|A|^2A\\
-2ic_1[\psi \p_y\tilde{A}+\phi \p_y a]+i\beta[a\p_y\tilde{\phi}+ A\p_y\psi],
\end{multline}
\begin{multline}\label{sm2}
\p_t \psi =D_1\p_x^2\psi+D_2\p_y^2 \psi -h\psi+\psi|\tilde{A}|^2+\phi(|\tilde{A}|^2-| A|^2)\\
 -c_{2}\mbox{Im} \left[ a^{\ast}\partial_{y}\tilde{A}+A^{\ast}\partial_{y} a\right]
\end{multline}
under the homogeneous Dirichlet boundary conditions and the initial conditions
$$
a(x,0)=a_0(x):=\tilde{A}_{0}(x)-A_{0}(x),\quad \psi_0(x)=\tilde{\phi}_0(x)-\phi_0(x),\quad x\in \Om.
$$
Now, we multiply \eqref{sm1} by $-\Delta a^{\ast}$ and \eqref{sm2} by $-\LL \psi$, add obtained relations and use the Young inequality

\begin{multline}\label{sm3}
\frac{1}{2}\frac{d}{dt}[\tau\lVert\nabla a\rVert^{2}+(\mathcal{L}\psi,\psi)]-\lVert\nabla a\rVert^{2}+\frac{1}{8}\lVert\Delta a \rVert^{2}+\frac{1}{2}\lVert\mathcal{L}\psi\rVert^{2}+h(\mathcal{L}\psi,\psi)\\
\leq 2(\phi^{4},\lvert a\rvert^{2})+2(\psi^{2}\lvert\tilde{\phi}+\phi\rvert^{2},\lvert\tilde{A}\rvert^{2})+2(\lvert\tilde{A}\rvert^{4},\lvert a\rvert^{2})\\
+2(\lvert a\rvert^{2}(\lvert\tilde{A}\rvert+\lvert A\rvert)^{2},\lvert A\rvert^{2})+8c_{1}^{2}(\lvert\psi\rvert^{2},\lvert\partial_{y}A\rvert^{2})+8c_{1}^{2}(\lvert\phi\rvert^{2},\lvert\partial_{y}A\rvert^{2})\\
+4\beta^{2}(\lvert a\rvert^{2},\lvert\partial_{y}\tilde{\phi}\rvert^{2})+4\beta^{2}(\lvert A\rvert^{2},\lvert\partial_{y}\psi\rvert^{2})+2(\psi^{2},\lvert A\rvert^{4})\\
+2(\phi^{2}\lvert a\rvert^{2},(\lvert\tilde{A}\rvert+\lvert A\rvert)^{2})+2c_{2}^{2}(\lvert a\rvert^{2},\lvert\partial_{y}\tilde{A}\rvert^{2})+2c_{2}^{2}(\lvert A\rvert^{2},\lvert\partial_{y}a\rvert^{2}).
\end{multline}
Next, we estimate terms on the right hand side of \eqref{sm3} by employing Young's inequality, the Ladyzhenskaya  and the Gagliardo-Nirenberg inequalities
\begin{multline}\label{E1}
2(\phi^{4},\lvert a\rvert^{2})\leq2\lVert\phi\rVert_{L^{4}}^{4}\lVert a\rVert_{L^{\infty}}^{2}\leq C_{1}\lVert a\rVert\lVert\Delta a\rVert\lVert\nabla\phi\rVert^{4}\\
\leq\varepsilon_{1}\lVert\Delta a\rVert^{2}+C(\varepsilon_{1})\lVert a\rVert^{2}\lVert\nabla \phi\rVert^{8},
\end{multline}
\begin{multline}\label{E2}
2(\psi^{2}\lvert\tilde{\phi}+\phi\rvert^{2},\lvert\tilde{A}\rvert^{2})\leq 4\lVert\phi\rVert_{L^{\infty}}^{2}(\lVert\tilde{\phi}\rVert_{L^{4}}^{2}+\lVert\phi\rVert^{2}_{L^{4}})\lVert\tilde{A}\rVert_{L^{4}}^{2}\\
\leq C_{2}\lVert\psi\rVert\lVert\Delta\psi\rVert(\lVert\nabla\tilde{\phi}\rVert^{2}+\lVert\nabla\phi\rVert^{2})\lVert\nabla\tilde{A}\rVert^{2}\\
\leq\varepsilon_{2}\lVert\mathcal{L}\psi\rVert^{2}+C(\varepsilon_{2})(\lVert\nabla\tilde{\phi}\rVert^{2}+\lVert\nabla\phi\rVert^{2})^{2}\lVert\nabla\tilde{A}\rVert^{4}\lVert\psi\rVert^{2},
\end{multline}
\begin{equation}\label{E3}
2(\lvert\tilde{A}\rvert^{4},\lvert a\rvert^{2})\leq\varepsilon_{1}\lVert\Delta a\rVert^{2}+C(\varepsilon_{1})\lVert a\rVert^{2}\lVert\nabla\tilde{A}\rVert^{8},\quad\quad\quad\quad\quad\quad\quad\quad\quad\quad
\end{equation}
\begin{multline}\label{E4}
2(\lvert a\rvert^{2}(\lvert\tilde{A}\rvert+\lvert A\rvert)^{2},\lvert A\rvert^{2})\\
\leq\varepsilon_{1}\lVert\Delta a\rVert^{2}+C(\varepsilon_{1})(\lVert\nabla\tilde{A}\rVert^{2}+\lVert\nabla A\rVert^{2})^{2}\lVert\nabla A\rVert^{4}\lVert a\rVert^{4},
\end{multline}
\begin{multline}\label{E5}
8c_{1}^{2}(\lvert\psi\rvert^{2},\lvert\partial_{y}A\rvert^{2})\leq C_{3}\lVert\psi\rVert\lVert\Delta\psi\rVert\lVert\nabla A\rVert^{2}\\
\leq\varepsilon_{2}\lVert\mathcal{L}\psi\rVert^{2}+C(\varepsilon_{2})\lVert\psi\rVert^{2}\lVert \nabla A\rVert^{4}.
\end{multline}
The rest of the terms on the right hand side of \eqref{sm3} can be estimated similarly.
Utilizing the estimates \eqref{E1}-\eqref{E5} in \eqref{sm3} and choosing $9\eb_1=\frac1{16}, \ 4\eb_2=\frac{1}{4}$ we arrive at the inequality
\begin{equation}\label{sm2a}
\frac d{dt}\left[\tau\|\Nx a\|^{2}+(\LL \psi,\psi)\right]+ \tau\|\Nx a\|^{2}+h(\LL \psi,\psi)\le M_{1}\left[\|a\|^2+\|\psi\|^2\right],
\end{equation}
where $M_1$ depends only on $\|\Nx A_0\|, \|\Nx \phi_0\|, \|\Nx \t A_0\|, \|\Nx \t \phi_0\|$ and $t$. The desired smoothing estimate \eqref{sma1} follows from the last inequality and the Lipschitz continuity of the semigroup  with respect to $[A,\phi]$ in $V^0$ established in \cite[Theorem II.1]{KKV}. 

It remains to show that $S(t)$ is H{\"o}lder continuous with respect to $t$, that it satisfies the condition b). Multiplying \eqref{a0} and \eqref{a1} by $\partial_{t}A^{\ast}$ and $\partial_{t}\phi$ respectively, and then using the estimates \eqref{chD0}, \eqref{A12} we deduce that
$$
\int_0^T\left[\|\p_t A(t)\|^2+\|\p_t \phi(t)\|^2\right]dt\le C_{T},\ \forall T>0,\ \forall [A,\phi]\in B_0,
$$
where $C_{T}$ depends only on $T$ and $B_0$. Thus for each $t_1,t_2\!\in\! [0,T]$ and $[A,\phi]\in\!B_0$
$$
[A(t_2)-A(t_1), \phi(t_2)-\phi(t_1)]=\int_{t_1}^{t_2}[\p_t A(t),\p_t\phi(t)]dt,
$$
and thanks to the Cauchy-Schwarz inequality we have
\begin{multline*}
\|[A(t_2)-A(t_1),\phi(t_2)-\phi(t_1)]\|_{V^0}=\int_{t_1}^{t_2}\|[\p_t A(t),\p_t\phi(t)]\|dt\\
\le |t_1-t_2|^{\frac12}\int_0^T\left[\|\p_t A(t)\|^2+\|\p_t \phi(t)\|^2\right]\le C_{T}|t_1-t_2|^{\frac12}.
\end{multline*}
Hence all conditions of the Theorem \ref{MZ} are satisfied. So we proved the following theorem:
\begin{theorem} If the conditions of the Theorem \ref{exun} are satisfied, then the semigroup $S(t): V^0\rw V^0$, $t\ge0$ possesses an exponential attractor.
\end{theorem}


\section{Feedback stabilization by finitely many Fourier modes}\label{sec:fs}
This section is devoted to the problem of global stabilization of the zero steady state of the system \eqref{a0}-\eqref{a1} by finitely many Fourier modes.
There are a number of papers on the problem of  stabilization of solutions of nonlinear PDE's  using finite-dimensional feedback controllers (see, e.g. \cite{AzTi}, \cite{BaTr}, \cite{BaWa}, \cite{Ch},  \cite{GuKa}, \cite{KaTi}, \cite{KaOz}, \cite{LuTi} and references therein). Our study of this question  is mainly inspired by  the results obtained in \cite{AzTi}, where the authors used various types of finite-dimensional controllers to stabilize the zero steady state solution  to semilinear parabolic equations. Accordingly, we study the following feedback control system
\begin{eqnarray}
\label{chp1}\tau\p_t A&=&A+\Delta A -\phi^2A-|A|^2A-2i c_1\phi\partial_y A\\
&&+i\beta A\partial_y \phi-\mu\sum\limits_{k=1}^N(A,w_k)w_k,\nonumber\\
\label{chp2}\partial_t \phi&=&D_1\partial^2_x\phi+D_2\partial^{2}_y\phi -h\phi +\phi\lvert A\rvert^{2}-c_{2}\mbox{Im}
\left[ A^{\ast}\partial_y A\right],
\end{eqnarray}
\be\label{chp3}
 A\big|_{\p \Om}=\phi\big|_{\p \Om}=0, \ \
A\big|_{t=0}=A_0, \ \phi\big|_{t=0}=\phi_0,
\ee
where $A_{0}, \phi_{0}\in L^{2}(\Omega)$ are given functions and the conditions on $\mu>0$ will be determined later. The proofs of the global existence and uniqueness for \eqref{chp1}-\eqref{chp3} and of the  estimate 
\begin{multline}\label{chpa1}
\frac12 \frac d{dt} \left[\tau\| A \|^2 +\delta_{0}\|\phi\|^2\right]-\|A\|^2+\|A\|^4_{L^4}\\
+
\delta_{0}D_0\| \nabla\phi\|^2+\delta_{0} h\|\phi\|^2 +\delta_0\| \nabla A\|^2
\le -\mu\sum\limits_{k=1}^N|(A,w_k)|^2,
\end{multline}
where $\delta_{0}$ and $D_{0}$ are as defined in \eqref{par1}-\eqref{par2}, are essentially the same with  the proof of Theorem \ref{exun} in Section \ref{sec:pr} and the estimate (2.9) in \cite{KKV}, respectively, which we will therefore omit.

Employing the inequality \eqref{PFN}, we infer from \eqref{chpa1} that
\begin{multline}\label{chpa2}
\frac12 \frac d{dt} \left[\tau\| A \|^2 +\delta_{0}\|\phi\|^2\right]+(\mu-1)\sum\limits_{k=1}^N|(A,w_k)|^2+\|A\|^4_{L^4}\\
+
\delta_{0}D_0\| \nabla\phi\|^2+\delta_{0} h\|\phi\|^2 +(\delta_0-\la_{N+1}^{-1})\| \nabla A\|^2
\le 0.
\end{multline}
Choosing $\mu \ge 1$ and $N$ large enough such that $\la_{N+1}^{-1}<\delta_{0}$ and using the Poincar{\'e} inequality, it follows from \eqref{chpa2} that
\be\label{A4}
\frac12 \frac d{dt} \left[\tau\| A \|^2 +\delta_{0}\|\phi\|^2\right]
+
\delta_{0}D_0\| \nabla\phi\|^2+\delta_{0} h\|\phi\|^2 +\delta_1\| A\|^2+\|A\|^4_{L^4}\le 0,
\ee
where $\delta_1:=\la_1(\delta_0-\la_{N+1}^{-1})$. From this inequality we infer that
$$
\frac12\frac d{dt} \left[\tau\| A \|^2 +\delta_{0}\|\phi\|^2\right]+m_0\left[\tau\| A \|^2 +\delta_{0}\|\phi\|^2\right]\le0,
$$
where $m_0=\min\{h+D_0\la_{1},\frac{\delta_{1}}{\tau}\}$. Thus we have
\be\label{est1}
\tau\| A(t) \|^2 +\delta_{0}\|\phi(t)\|^2\le e^{-2m_0t}\left[\tau\| A_0 \|^2 +\delta_{0}\|\phi_0\|^2\right],
\ee
which proves that all solutions of the controlled problem tend to the zero steady state as $t\rw \infty$ with an exponential rate in $L^{2}(\Omega)$ whenever
\be\label{exp1}
\mu \ge 1\quad \mbox{and}\quad \la_{N+1}^{-1}< \delta_0.
\ee
Finally, we would like to note that, the estimates obtained above allow us to claim the global unique solvability of the problem \eqref{chp1}-\eqref{chp3} in 
\begin{equation}
C(0,T; V^1)\cap L^2(0,T; V^2),\quad \forall T>0.\nonumber
\end{equation} 
Moreover the smoothing estimate \eqref{sm2a} is also valid for solutions of the problem \eqref{chp1}-\eqref{chp3}. Thus the following theorem holds true:

\begin{theorem}
If the number  $N$ and $\mu$ satisfy the conditions \eqref{exp1},
then the solution $[A(t),\phi(t)]$  of the problem \eqref{chp1}-\eqref{chp3} tends to zero stationary state  with an exponential rate in $V^{1}$.
\end{theorem}


\section{1D Chevron pattern equations}
In this section, we study 1D version of the system \eqref{a0}-\eqref{a2}, reduced to one spatial dimension by eliminating the dependence on $y$ and taking $x$ to lie in the interval $\Omega=[0,L]$.
\be\label{1Dch}
\begin{cases}
\tau\p_t A-\p_x^2A-A+|A|^2A +\phi^2A=0, \\
\p_t \phi -D_1\p_x^2 \phi +h\phi -|A|^2\phi =0,\quad \ x\in (0,L),\quad t>0, \\
 A\big|_{x=0}=A\big|_{x=L}=\phi\big|_{x=0}=\phi\big|_{x=L}=0, \\
A\big|_{t=0}=A_0,\quad \phi\big|_{t=0}=\phi_0,
\end{cases}
\ee
We take the scalar product of the first equation in \eqref{1Dch} with $A^{\ast}$, and the second equation with $\phi$. Then after a series of integration by parts, we add the resulting equations and obtain
 \be\label{1D1}
\frac12\frac d{dt}\left[\tau\|A\|^2\!+\|\phi\|^2\right]+\|\p_xA\|^2\!+D_{1}\|\p_x\phi\|^2\!-\|A\|^2+\|A\|_{L^4}^4\!+h\lVert\phi\rVert^{2}=0.
\ee
Employing the Poincar{\'e} inequality and the following inequality
$$
\|A\|^2\le \frac{L}{2}+\frac12\|A\|_{L^4((0,L))}^4,
$$
in \eqref{1D1}, we get
\be\label{1D2}
\frac d{dt} \left[\tau\|A(t)\|^2+\|\phi(t)\|^2\right]+\alpha\left[\tau\|A(t)\|^2+\| \phi(t)\|^2\right]+\|A\|_{L^4}^4+2h\lVert\phi\rVert^{2}\le L,\nonumber
\ee
where $\alpha=\min\{\frac{2\pi^{2}}{\tau L^{2}}+\frac{2}{\tau},\frac{2\pi^{2}D_{1}}{L^{2}}\}$. Then
\be\label{Dis1}
\tau\|A(t)\|^2+\|\phi(t)\|^2\le e^{-\alpha t}\left[\tau\|A_0\|^2+\|\phi_0\|^2\right]+\frac{L}{\alpha},
\ee
due to Gronwall lemma, and integrating the resulting inequality in $t$ over $(0,T)$, yields
\begin{equation}\label{Dis1a}
\int_0^T\!\left[\tau\|A(t)\|^2\!+\|\phi(t)\|^2\right]dt
\le\frac{1-e^{-\alpha T}}{\alpha}(\tau \|A_0\|^2+\|\phi_0\|^2)+\frac{LT}{\alpha}, 
\end{equation}
$\forall T>0$.
\begin{remark}
We can prove that the semigroup $S(t): V^{0}\rightarrow V^{0}$ generated by the problem \eqref{1Dch} possesses an exponential attractor, using a similar argument to the one in section \ref{sec:ea}.\\
It is easy to see that
\begin{equation}
\frac d{dt} \Lambda (t)=-2\|\p_tA(t)\|^2-2\|\p_t\phi(t)\|^2,\nonumber
\end{equation}
where
\begin{multline}
\Lambda(t):=\\
  \|\p_x A(t)\|^2+D_1\|\p_x \phi(t)\|^2+\frac12\| A(t)\|^4_{L^4}+ (\phi^2(t),|A(t)|^2)-
\|A(t)\|^2+h\|\phi(t)\|^2.\nonumber
\end{multline}
So $\Lambda(t)$ is a Lyapunov function for the system. Therefore
 the global attractor of the system consists of stationary states and trajectories joining them (if  $L\gg 1$),  (see, e.g. \cite{Lad}). 
 
Let us note that the system \eqref{1Dch} has a rich family of stationary states, including stationary states of the form $[W,0]$ (when $L\gg 1$), where $W$ is a stationary state of 1D Ginzburg-Landau equation:
\begin{displaymath}
\left\{ \begin{array}{l}
-W''-W+\lvert W\rvert^{2}W=0,\quad x\in(0,L),\nonumber\\
W(0)=W(L)=0,\nonumber
\end{array} \right.
\end{displaymath}
(see, e.g., \cite{BaVi}).
 
 Moreover  the system possesses an inertial manifold (see, e.g., \cite{KoZe}, \cite{KoZe2} and references therein).

\end{remark}

\subsection{Feedback stabilization}
The zero solution of the system \eqref{1Dch} is unstable for $L\gg1$, which can easily be seen by setting $\phi\equiv0$ and $A=O(\epsilon)$, where $\epsilon\ll 1$. In this case, the solution of the  linearization of the first equation in \eqref{1Dch} around zero 
\begin{equation}
\tau \partial_{t}A-\partial_{x}^{2}A-A=0,\nonumber
\end{equation}
has the form
\begin{equation}
A=\sum_{k=0}^{\infty}\alpha_{k}\sin\left(\frac{k\pi}{L}x\right),\nonumber
\end{equation}
where
\begin{equation}
\alpha_{k}(t)=\alpha_{k}(0)e^{-\left(\frac{k^{2}\pi^{2}-L^{2}}{\tau L^{2}}\right)t},\nonumber
\end{equation}
which implies that for $L\gg1$ as $t\rightarrow\infty$, $A$ does not go to zero. Motivated by this observation, we study the feedback stabilization of the 1D system by using the a priori estimates \eqref{Dis1} and \eqref{Dis1a}. Contrary to section \ref{sec:fs} where it is the zero steady state that is stabilized, in the 1D case we can stabilize any time-dependent, not necessarily stationary, solution.

We start with deriving the uniform estimates of solutions to \eqref{1Dch}. Multiplying the first equation in \eqref{1Dch} by $-\p_x^2 A^{\ast}$ in $L^2$, and utilizing the inequality
$$
-Re (|A|^2A, \p_x^2A)\ge (|A|^2,|\p_xA|^2),
$$
gives us
\begin{multline}\label{best4}
\frac{\tau}{2}\frac d{dt}\|\p_x A\|^2+\|\p^2_x A\|^2-\|\p_x A\|^2+(|A|^2,|\p_x A|^2)+(\phi^2,|\p_x A|^2)\\
\le 2|(\phi \p_x\phi, A^{\ast}\p_xA)|\le \eb_1(|A|^2,|\p_xA|^2)+C_{\eb_1}(|\phi|^2,|\p_x\phi|^2).
\end{multline}
We estimate the last terms on the right hand side of \eqref{best4} and the term $\lVert\partial_{x}A\rVert^{2}$, by using the Gagliardo-Nirenberg inequality \eqref{GN2}
$$
\|u\|_{L^4}\le C\|u\|^{7/8}\|\p^2_x u\|^{1/8},\quad \|\p_x u\|_{L^4}\le C \|u\|^{3/8}\|\p^2_x u\|^{5/8}
$$
and Young's inequality, as follows
\begin{eqnarray}
\|\p_xA\|^2&\le& \frac1{4\eb_1}\|A\|^2+ \eb_{1} \|\p^2_xA\|^2,\nonumber\\
C_{\varepsilon_{1}}(|\phi|^2,|\p_x\phi|^2)&\le& C_{\varepsilon_{1}}\|\phi\|^2_{L^4}\|\p_x \phi\|_{L^4}^2\le \label{Fes1}C_{\varepsilon_{1}} C^{2}\|\phi\|^{5/2}\|\p_x^2\phi\|^{3/2}\\
&\le&\eb_2\|\p_x^2\phi\|^2+C_{\eb_2}\|\phi\|^{10}.\nonumber
\end{eqnarray}
So, we have
\begin{multline}\label{best5}
\frac\tau2\frac d{dt}\|\p_x A\|^2+(1-\eb_1)\|\p^2_x A\|^2+(1-\eb_1)(|A|^2,|\p_x A|^2)\\
\le \eb_2\|\p_x^2\phi\|^2+C_{\varepsilon_{2}}\|\phi\|^{10}+\frac1{4\eb_1}\|A\|^2.
\end{multline}
Next we multiply the second equation in \eqref{1Dch} by $-\p_x^2\phi$ in $L^2$:
\begin{multline}\label{best5a}
\frac12\frac d{dt}\|\p_x \phi\|^2+D_1\|\p^2_x \phi\|^2+h\|\p_x\phi\|^2
\le
(|A|^2,|\p_x\phi|^2)\\
+2(|A||\p_xA|,|\phi||\p_x\phi|)
\le 2(|A|^2,|\p_x\phi|^2)+(\phi^2,|\p_x A|^2)
\end{multline}
We estimate the first term on the right hand side of \eqref{best5a}, again by using Gagliardo-Nirenberg inequality
\begin{multline}\label{best5b}
2(|A|^2,|\p_x\phi|^2)\le 2\|A\|^2_{L^4}\|\p_x \phi\|_{L^4}^2\le 2C^{4}\|A\|^{\frac{7}{4}}\|\p^2_xA\|^{\frac{1}{4}}
\|\phi\|^{\frac{3}{4}}\|\p^2_x\phi\|^{\frac{5}{4}}\\
\le \eb_{1}\|\p^2_xA\|^2
+\eb_{2}\|\p^2_x\phi\|^2+C_{2}(\varepsilon_{1})\|A\|^{14}+
C_{2}(\varepsilon_{2})\|\phi\|^{6}.
\end{multline}
Combining \eqref{Fes1} and \eqref{best5b} with \eqref{best5a}, it follows that
\begin{multline*}
\frac12\frac d{dt}\|\p_x \phi\|^2+(D_1-2\eb_2)\|\p^2_x \phi\|^2+h\|\p_x\phi\|^2
\\
\le \eb_{1}\|\p^2_xA\|^2+C_{\varepsilon_{2}}\lVert\phi\rVert^{10}+C_{2}(\varepsilon_{1})\lVert A\rVert^{14}+C_{2}(\varepsilon_{2})\lVert\phi\rVert^{6}.
\end{multline*}
Adding the above inequality to \eqref{best5} we get
\begin{multline}\label{best5c}
\frac12\frac d{dt}\left[\tau \|\p_x A\|^2+\|\p_x \phi\|^2\right]+(1-2\eb_1)\|\p^2_x A\|^2+(1-\eb_1)(|A|^2,|\p_x A|^2)
\\+(D_1-3\eb_2)\|\p^2_x \phi\|^2+h\|\p_x\phi\|^2\le 2C_{\varepsilon_{2}}\|\phi\|^{10}\\
+\frac1{4\eb_1}\|A\|^2+C_2(\eb_1)\|A\|^{14}+
C_2(\eb_2)\|\phi\|^{6}\nonumber
\end{multline}

The resulting inequality after setting $\eb_1=\frac14, \eb_2=\frac14D_1$ in the above inequality implies the dissipativity of the system
in the phase space $V^{1}$. More precisely the following inequality holds true
\begin{equation}\label{Dis2}
\frac d{dt}\E_1(t)+d_0\E_1(t)\le Q_1(\|A\|) + Q_2(\|\phi)\|),
\end{equation}
where $\E_1(t):=\tau\|\p_x A\|^2+\|\p_x \phi\|^2$, $d_{0}=\frac{\pi^{2}}{L^{2}}\min\{\frac{1}{2\tau},\frac{D_{1}}{4}\}$ and
$Q_1(\cdot)$, $Q_2(\cdot)$ are monotone functions.

Assume that $[A,\phi]$ is an arbitrary given and possibly time-dependent  solution of the problem \eqref{1Dch} and consider the following feedback control system

\be\label{con1}
\begin{cases}
\tau\p_t \t A-\p_x^2\t A-\tilde{A}+|\t A|^2\t A +{\t \phi}^2\t A=-\mu \sum\limits_{k=1}^N(\t A,w_k)w_k(x), \\
\p_t \t\phi -D_1\p_x^2 \t \phi +h\t\phi -|\t A|^2\t \phi =0,\quad x\in (0,L),\quad t>0, \\
\t A\big|_{x=0}=\t A\big|_{x=L}=\t\phi\big|_{x=0}=\t\phi\big|_{x=L}=0, \\
\t A\big|_{t=0}=\t A_0, \ \t\phi\big|_{t=0}=\t\phi_0,
\end{cases}
\ee
where $w_k(x)=\sin (\frac{k\pi}{L} x), k=1,2,...$. Then the pair of functions $[a,\psi]:=[\t A-A,\t \phi-\phi]$ is a solution of the problem
\be\label{con2}
\begin{cases}
\tau\p_t a-\p_x^2a-a+|\t A|^2\t A-| A|^2A+\phi^2 a+({\t\phi}^2-\phi^2)\t A\\
=-\mu_{1} \sum\limits_{k=1}^{N_{1}}(a,w_k)w_k(x), \\
\p_t \psi -D_1\p_x^2 \psi +h\psi -| A|^2\psi =(|\t A|^2-|A|^2)\t \phi
\\-\mu_{2} \sum\limits_{k=1}^{N_{2}}(\psi,w_k)w_k(x),\quad x\in (0,L),\quad t>0, \\
a\big|_{x=0}=a\big|_{x=L}=\psi\big|_{x=0}=\psi\big|_{x=L}=0, \\
\t a\big|_{t=0}=a_0, \ \psi\big|_{t=0}=\psi_0,
\end{cases}
\ee
where $a_0=\t{A}_0-A_0, \ \psi_0=\t\phi_0- \phi_0.$  Taking the inner product of the first equation in \eqref{con2} with $a^{\ast}$ and of the second equation with $\psi$, we get
\begin{multline}\label{est:e1}
\frac\tau2\frac d{dt} \|a\|^2+\|\p_xa\|^2-\|a\|^2+(\phi^2,|a|^2)+(\lvert\tilde{A}\rvert^{2},\lvert a\rvert^{2})\\
+Re ({\t\phi}^2-\phi^2,\tilde{A}a)\le (\lvert\tilde{A}\rvert\lvert A\rvert,\lvert a\rvert^{2})+(\lvert A\rvert^{2},\lvert a\rvert^{2})-\mu_1\sum_{k=1}^{N_1}|(a,w_k)|^2,
\end{multline}

\begin{multline}\label{est:e2}
\frac12\frac d{dt} \|\psi\|^2+D_1\|\p_x\psi\|^2+h\|\psi\|^2=\\
(|A|^2,|\psi|^2)+(|A|^2-|\t A|^2,\t\phi\psi)-\mu_2 \sum\limits_{k=1}^{N_2}(\psi,w_k)^2.
\end{multline}
We use the Sobolev inequality $\lVert u\rVert_{L^{\infty}}\leq \lVert \partial_{x}u\rVert$, H\"older's inequality and Young's inequality \eqref{Young} to produce the following estimates
\begin{multline}\label{est:e3}
\lvert(\phi^{2}-\tilde{\phi}^{2},\tilde{A}a)\rvert\leq(\lvert\psi\rvert\lvert\phi+\tilde{\phi}\rvert,\lvert\tilde{A}\rvert\lvert a\rvert)\\
\leq \frac{1}{2}(\lvert a\rvert^{2},\lvert A\rvert^{2})+\frac{1}{2}(\lvert\psi\rvert^{2},\lvert\phi+\tilde{\phi}\rvert^{2})\\
\quad\quad\quad\quad\quad\quad\quad\leq\frac{1}{2}\lVert A\rVert_{L^{\infty}}^{2}\lVert a\rVert^{2}+(\lVert\phi\rVert_{L^{\infty}}^{2}+\lVert\tilde{\phi}\rVert_{L^{\infty}}^{2})\lVert\psi\rVert^{2}\\
\leq\frac{1}{2}\lVert\partial_{x}A\rVert^{2}\lVert a\rVert^{2}+(\lVert\partial_{x}\phi\rVert^{2}+\lVert\partial_{x}\tilde{\phi}\rVert^{2})\lVert\psi\rVert^{2},
\end{multline}
\begin{eqnarray}
\label{est:e3}&&(\lvert A\rvert^{2},\lvert a\rvert^{2})\!\leq\!\lVert\partial_{x}A\rVert^{2}\lVert a\rVert^{2},\quad (\lvert \tilde{A}\rvert\lvert A\rvert,\lvert a\rvert^{2})\!\leq\!(\lVert\partial_{x}A\rVert^{2}+\lVert\partial_{x}\tilde{A}\rVert^{2})\lVert a\rVert^{2},\\
\label{est:e4}&&(\lvert A\rvert^{2},\psi^{2})\leq\lVert\partial_{x}A\rVert^{2}\lVert\psi\rVert^{2},
\end{eqnarray}
\begin{eqnarray}\label{est:e7}
\quad\quad\lvert(\lvert A\rvert^{2}-\lvert\tilde{A}\rvert^{2},\tilde{\phi}\psi)\rvert&\leq&(\lvert a\rvert\lvert A+\tilde{A}\rvert,\lvert\tilde{\phi}\rvert\lvert\psi\rvert)\\
&\leq&\frac{1}{2}(\lvert A+\tilde{A}\rvert^{2},\lvert a\rvert^{2})+\frac{1}{2}(\lvert\tilde{\phi}\rvert^{2},\lvert\psi\rvert^{2})\nonumber\\
&\leq&(\lVert\partial_{x}A\rVert^{2}+\lVert\partial_{x}\tilde{A}\rVert^{2})\lVert a\rVert^{2}+\frac{1}{2}\lVert\partial_{x}\tilde{\phi}\rVert^{2}\lVert\psi\rVert^{2}.\nonumber
\end{eqnarray}

Adding \eqref{est:e1} to \eqref{est:e2} and using the estimates \eqref{est:e3}-\eqref{est:e7}, we get
\begin{multline}\label{cona1}
\frac12\frac d{dt}\left[ \tau \|a\|^2+\|\psi\|^2\right]+\|\p_xa\|^2+D_1\|\p_x\psi\|^2-\|a\|^2 +h\|\psi\|^2\\
\le \left(\frac{7}{2}\|\p_x A\|^2+2\|\p_x \t A\|^2 \right)\|a\|^2+\left(\|\p_x A\|^2+\|\p_x \phi\|^2 +\frac{3}{2}\|\p_x \t \phi\|^2\right)\|\psi\|^2\\
-\mu_1\sum_{k=1}^{N_1}|(a,w_k)|^2-\mu_2\sum_{k=1}^{N_2}|(\psi,w_k)|^2.
\end{multline}
Since $$
\|\p_x \t A\|^2, \|\p_x A\|^2, \|\p_x \phi\|^2, \|\p_x \t \psi\|^2\le M_0,
$$
we have
\begin{multline*}
\frac12\frac d{dt}\left[ \tau \|a\|^2+\|\psi\|^2\right]+\|\p_xa\|^2+D_1\|\p_x\psi\|^2-\|a\|^2 +h\|\psi\|^2\\
\le 6M_0(\|a\|^2+\|\psi\|^2)
-\mu_1\sum_{k=1}^{N_1}|(a,w_k)|^2-\mu_2\sum_{k=1}^{N_2}|(\psi,w_k)|^2.
\end{multline*}
We rewrite the above inequality in the following form
\begin{multline}\label{con2a}
\frac12\frac d{dt}\left[ \tau \|a\|^2+\|\psi\|^2\right]+\|\p_xa\|^2+D_1\|\p_x\psi\|^2\\
+(\mu_1-1-6M_0)\sum_{k=1}^{N_1}|(a,w_k)|^2
+(\mu_2+h-6M_0)\sum_{k=1}^{N_2}|(\psi,w_k)|^2\\
\leq (1+6M_{0})\sum_{k=N_{1}+1}^{\infty}\lvert(a,w_{k})\rvert^{2}+6M_{0}\sum_{k=N_{2}+1}^{\infty}\lvert(\psi,w_{k})\rvert^{2}\\
\leq (1+6M_{0})\lambda_{N_{1}+1}^{-1}\lVert\partial_{x}a\rVert^{2}+6M_{0}\lambda_{N_{2}+1}^{-1}\lVert\partial_{x}\psi\rVert^{2}.
\end{multline}
Assuming that
\begin{equation}\label{con:star2}
(1+6M_{0})\lambda_{N_{1}+1}^{-1}\leq\frac{1}{2},\ 6M_{0}\lambda_{N_{2}+1}^{-1}\leq\frac{D_{1}}{2},\ \mu_{1}\geq 1+6M_{0},\ \mu_{2}+h\geq 6M_{0}
\end{equation}
and using the inequality \eqref{PFN} in \eqref{con2}, we derive the desired inequality
$$
\frac d{dt}\left[ \tau \|a\|^2+\|\psi\|^2\right]+\la_1\left[ \|a\|^2+D_{1}\|\psi\|^2\right]\leq 0
$$
which implies that
$$
\tau \|a(t)\|^2+\|\psi(t)\|^2\le \left[\tau \|a_0\|^2+\|\psi_0\|^2\right]e^{-\la_1r_0t},
$$
where $r_0:=\min\{\frac1\tau,D_{1}\}$.

Hence we proved the following
\begin{theorem}\label{1Ds} Let $[A(t),\phi(t)]$ be a given solution of the problem \eqref{1Dch}.
If $\mu_{1}$, $\mu_{2}$, $N_{1}$ and $N_{2}$ are so large that \eqref{con:star2} is satisfied then each solution of the controlled problem \eqref{con1} tends to $[A(t),\phi(t)]$ as $t\rw \infty$ in $V^{0}$ with an exponential rate.
\end{theorem}
\begin{remark} Let us note that analog of the Theorem \ref{1Ds} is true also for the 2D system of equations. It can be proved in a similar way thanks to the estimates \eqref{A12}, \eqref{chD2}.
\end{remark}


\subsection{Numeric results}

We consider the following one-dimensional version of the Chevron problem \eqref{chp1}--\eqref{chp2} over the interval $0\leq x\leq L$, where the dependence on the $y$-direction has been eliminated:
\begin{align}
&\tau\p_t A-\p_x^2A-A+|A|^2A +\phi^2A=0, \\
&\p_t \phi -D_1\p_x^2 \phi +h\phi -|A|^2\phi =0,\quad x\in (0,L),\quad t>0, \\
 &A\big|_{x=0}=A\big|_{x=L}=\phi\big|_{x=0}=\phi\big|_{x=L}=0, \\
&A\big|_{t=0}=A_0, \ \phi\big|_{t=0}=\phi_0,
\end{align}
In order to study the problem numerically, we introduce the following semi-implicit finite-differences scheme with uniform space and time discretization steps $\delta x = L/N$ and $\delta t>0$ so that $A^k_i \approx A(i\delta x,k \delta t)$, for $0\leq i\leq N$:
\begin{align}\label{eq:schemeCb}
&\tau \frac{A_i^{k+1}-A_i^k}{\delta t} - \frac{A_{i-1}^{k+1} - 2A_i^{k+1} + A_{i+1}^{k+1}}{\delta x^2} + (H^k_+)_iA_i^{k+1}= (H^k_-)_i A_i^k,\\
&\frac{\phi_i^{k+1}-\phi_i^k}{\delta t} - D_1\frac{\phi_{i-1}^{k+1} - 2\phi_i^{k+1} + \phi_{i+1}^{k+1}}{\delta x^2} +  (G^k_+)_i \phi_i^{k+1}= (G^k_-)_i \phi_i^k,\\
 & (H^k_+)_i = \lvert A^k_i\rvert^2 + \lvert \phi^k_i\rvert^2,  \quad(H^k_-)_i = 1,\\
 & (G^k_+)_i = h, \quad(G^k_-)_i = \lvert A^k_i\rvert^2,\\
& A^{k+1}_0 = A^{k+1}_N = \phi^{k+1}_0 = \phi^{k+1}_N = 0,\\\
\label{eq:schemeCe}&A^0_i = A_0(i \delta x), \phi^0_i = \phi_0(i \delta x), 0\leq i \leq N\,.
\end{align}
This leads to a linear system of the form:
\begin{align}
\label{eq:linsysA}&\frac{\tau}{\delta t}(A^{k+1}_h - A^k_h) + \Delta_h A^{k+1}_h + \mathbf{H}^k_+ A^{k+1}_h = \mathbf{H}^k_- A^k_h,\\
\label{eq:linsysP}&\frac{1}{\delta t}(\phi^{k+1}_h - \phi^k_h) + D_1 \Delta_h \phi^{k+1}_h + \mathbf{G}^k_+ \phi^{k+1}_h = \mathbf{G}^k_- \phi^k_h,
\end{align}
where the subscript $h$ denotes a vector of the form $f_h = \{ f_i\}_{i=0}^{N}$. The matrix $\Delta_h$ is tridiagonal and positive definite, with diagonal entries $(\Delta_h)_{ii} = \frac{2}{\delta x^2}$ and off-diagonal entries $(\Delta_h)_{i,i+1} = (\Delta_h)_{i-1,i} = -\frac{1}{\delta x^2}$, reflecting the Dirichlet boundary conditions at $x=0$ and $x=L$. The matrices $\mathbf{H}^k_\pm$ and $\mathbf{G}^k_{\pm}$ are diagonal with the non-negative vectors $H^k_\pm$ and $G^k_\pm$ resp.~in the main diagonal, and therefore symmetric and positive semi-definite. By construction, $H^k_+ - H^k_- = H^k_h = -1 + \lvert A^k_h\rvert^2 + \lvert \phi^k_h\rvert^2$ and $G^k_+ - G^k_- = G^k_h = h - \lvert \phi^k_h\rvert^2$.

Although a full analysis is beyond the scope of this paper, the behavior of the proposed scheme is captured by the following estimates for the growth of the $L^2$ vector norm $\lVert u_h\rVert_2^2 = \sum_{i=0}^N \lvert u_i\rvert^2$ and the discrete $H^1$ semi-norm $u_h^\dagger \Delta_h u_h$ of the solution:
\begin{multline}
  \frac{1}{2}\left(1-\frac{\delta t}{\tau}\lVert H^k_-\rVert_\infty\right)\lVert A^{k+1}_h\rVert^2  +\frac{\delta t}{\tau} (A^{k+1}_h)^\dagger \Delta_h A^{k+1}_h  \\
+\frac{\delta t}{\tau} (A^{k+1}_h)^\dagger\mathbf{H}_+^k A^{k+1}_h \leq \frac{1}{2}\left(1+ \frac{\delta t}{\tau}\lVert H^k_-\rVert_\infty\right)\lVert A^{k}_h\rVert^2,
\end{multline}
\begin{multline}
  \frac{1}{2}\left(1 - \delta t\lVert G^k_-\rVert_\infty\right)\lVert \phi^{k+1}_h\rVert^2  + \delta t D_1 (\phi^{k+1}_h)^\dagger \Delta_h \phi^{k+1}_h  \\
+\delta t (\phi^{k+1}_h)^T \mathbf{G}_+^k \phi^{k+1}_h \leq \frac{1}{2}\left(1+ \delta t \lVert G^k_-\rVert_\infty\right)\lVert \phi^{k}_h\rVert^2 \,.
\end{multline}
and the energy-type inequalities:
\begin{multline}
\tau \delta t\lVert \tfrac{A^{k+1}_h - A^k_h}{\delta t}\rVert_2^2 + \frac{1}{2}(A^{k+1}_h)^\dagger \Delta_h A^{k+1}_h +  \frac{1}{2} (A^{k+1}_h)^\dagger \mathbf{H}^{k+1}_h A^{k+1}_h \leq\\
   \frac{\delta t}{2}\lVert A^{k+1}_h\rVert^2 \lVert\tfrac{H^{k+1}_h - H^k_h}{\delta t}\rVert_\infty +  \frac{1}{2}(A^k_h)^\dagger \Delta_h A^{k}_h +   \frac{1}{2} (A^{k}_h)^\dagger \mathbf{H}^k_h A^{k}_h
\end{multline}
\begin{multline}
\delta t\lVert \tfrac{\phi^{k+1}_h - \phi^k_h}{\delta t}\rVert_2^2 + \frac{D_1}{2}(\phi^{k+1}_h)^T \Delta_h \phi^{k+1}_h +  \frac{1}{2} (\phi^{k+1}_h)^T \mathbf{G}^{k+1}_h \phi^{k+1}_h \leq\\
   \frac{\delta t}{2}\lVert \phi^{k+1}_h\rVert^2 \lVert\tfrac{G^{k+1}_h - G^k_h}{\delta t}\rVert_\infty +  \frac{D_1}{2}(A^k_h)^T \Delta_h \phi^{k}_h +   \frac{1}{2} (\phi^{k}_h)^T \mathbf{G}^k_h \phi^{k}_h
\end{multline}
The first two estimates follow from multiplying the equations \eqref{eq:linsysA}-\eqref{eq:linsysP} from the left with $(A^{k+1}_h)^\dagger$ and $(\phi^{k+1}_h)^T$ respectively. Likewise the other two estimates follow from left multiplication with $(A^{k+1}_h - A^k_h)^\dagger$ and $(\phi^{k+1}_h-\phi^k_h)^T$. Note also the use of the infinity vector norm $\lVert u_h\rVert_\infty = \max_i \lvert u_i\rvert$.

In matrix form, the system \eqref{eq:linsysA}-\eqref{eq:linsysP} can be written as
\begin{align}
\label{eq:discrete_systemA}&\mathbf{M}_A^k A^{k+1}_h = \mathbf{L}_A^k A^k_h, \\
\label{eq:discrete_systemP}&\mathbf{M}_\phi^k \phi^{k+1}_h = \mathbf{L}_\phi^k \phi_h^k,
\end{align}
with
\begin{align}
 & \mathbf{M}_A^k = I + \frac{\delta t}{\tau}\Delta_h + \frac{\delta t}{\tau} \mathbf{H}^k_+, &L_A^k =  I  + \frac{\delta t}{\tau} \mathbf{H}^k_-\\
& \mathbf{M}_\phi^k = I + \delta t D_1 \Delta_h + \delta t \mathbf{G}^k_+, & L_\phi^k = I + \delta t \mathbf{G}^k_-\,.
\end{align}
The matrices $M_A^k$ and $M_\phi^k$ are symmetric, unconditionally positive-definite and sparse, and hence the system \eqref{eq:discrete_systemA}-\eqref{eq:discrete_systemP} can be solved efficiently, either directly via a suitable factorization routine, or iteratively via for instance the Conjugate Gradient method.

For the stabilization with $K$ discrete modes $W^j_h\in \mathbb{R}^N$, $1\leq j\leq K$ stacked in an $N\times K$ matrix $\mathbf{W}_{ij} = (W^j_h)_i$, we extend the scheme \eqref{eq:linsysA}-\eqref{eq:linsysP} to
\begin{align}\label{eq:stab}
&\frac{\tau}{\delta t}(A^{k+1}_h - A^k_h) + \Delta_h A^{k+1}_h + \mathbf{H}^k_+A^{k+1}_h =\mathbf{H}^k_-A^{k}_h -\mu \mathbf{W}\mathbf{W}^T A^{k+1}_h,\\
&\frac{1}{\delta t}(\phi^{k+1}_h - \phi^k_h) + D_1 \Delta_h \phi^{k+1}_h + \mathbf{G}^k_+ \phi^{k+1}_h =\mathbf{G}^k_- \phi^{k}_h,
\end{align}
and correspondingly, in matrix form,
\begin{align}\label{eq:discrete_stab}
&(\mathbf{M}_A^k + \tfrac{\mu\delta t}{\tau} \mathbf{W}\mathbf{W}^T) A^{k+1}_h = \mathbf{L}_A^k A^k_h, \\
&\mathbf{M}_\phi^k \phi^{k+1}_h = \mathbf{L}_\phi^k \phi_h^k \,.
\end{align}
The stabilized matrix $\mathbf{M}_{A,W}^k = \mathbf{M}_A^k +  \frac{\mu\delta t}{\tau}\mathbf{W}\mathbf{W}^T$ is not sparse, but can still be inverted via an iterative method, such as CG, as a low-rank update of a sparse matrix for the extra cost of $K$ vector-vector products per matrix-vector product.

Let $0<\lambda_1<\ldots<\lambda_N$ be the eigenvalues of $\Delta_h$ and $V^j_h$ corresponding eigenvectors, mutually orthogonal and normalized so that $\lVert V^j_h\rVert_2 =1$. Then
\begin{equation}
\lVert (\mathbf{M}^k_A)^{-1}\mathbf{L}^k_A V^j_h\rVert_2 \leq \frac{1+\tfrac{\delta t}{\tau}\lVert {H}^k_-\rVert_\infty}{1+\tfrac{\delta t}{\tau}\lambda_j}
\end{equation}
and likewise for $\phi$. This implies that modes with $\lambda_j>\lVert {H}^k_-\rVert_\infty$ are damped (independent of the time step $\delta t$ notably). For the stabilized system on the other hand,
\begin{equation}
\lVert (\mathbf{M}^k_{A,W})^{-1}\mathbf{L}^k_A V^j_h\rVert_2
\leq \frac{1+\tfrac{\delta t}{\tau}\lVert {H}^k_-\rVert_\infty}{1+\tfrac{\delta t}{\tau}\lambda_j + \tfrac{\mu\delta t}{\tau}\lVert \mathbf{W}^T V^j_h\rVert_2^2}
\end{equation}
since the eigenvectors $V^j_h$ constitute an orthonormal basis and threfore
\begin{multline*}
 \lVert \mathbf{M}^k_{A,W}V^j_h\rVert_2^2 = \sum_{l=1}^N \lvert (V^l_h)^T\mathbf{M}^k_{A,W} V^j_h\rvert^2\\
 \geq \lvert (V^j_h)^T\mathbf{M}^k_{A,W} V^j_h\rvert^2 = \left(1+\tfrac{\delta t}{\tau}\lambda_j + \tfrac{\mu\delta t}{\tau}\lVert \mathbf{W}^T V^j_h\rVert_2\right)^2,
\end{multline*}
It follows that we have full damping if $\lambda_j + \mu\lVert \mathbf{W}^T V^j_h\rVert_2^2> \lVert {H}^k_-\rVert_\infty$ for all modes with $\lambda_j\leq \lVert {H}^k_-\rVert_\infty$.

In the Chevron system \eqref{eq:schemeCb}-\eqref{eq:schemeCe}, we have $(H^k_-)_i=1$, for all $0\leq i\leq N$, and so $\lVert H^k_-\rVert_\infty =1$ at all times $t^k$. This implies that we can choose a fixed $\mu>0$ and $\mathbf{W}$ such that $\lambda_j + \mu\lVert \mathbf{W}^T V^j_h\rVert_2^2> 1$ for all $\lambda_j<=1$, and then $\lVert \mathbf{M}^k_{A,W}V^j_h\rVert_2<1$ for all $V^j_h$, and so  $\lVert A^{k+1}_h\rVert_2 = \lVert \mathbf{M}^k_{A,W}A^k_h\rVert_2 <  \lVert A^k_h\rVert_2$. On the other hand $\lVert G^k_-\rVert_\infty = \lVert \lvert A^k_h\rvert^2 \rVert_\infty \leq \lVert A^k_h\rVert_2^2 \leq \lVert A^0_h\rVert_2^2$ because of the stabilisation of $A$. So if in addition  $\lambda_j + \mu\lVert \mathbf{W}^T V^j_h\rVert_2^2> \lVert A^0_h\rVert_2^2$ for any $\lambda_j \leq \lVert A^k_0\rVert_2^2$, then $\lVert \phi^{k+1}_h\rVert_2 < \lVert \phi^{k}_h\rVert_2$. We conclude that: \begin{quote}\textbf{if $\mu,\epsilon>0$ and $\mathbf{W}\in\mathbb{R}^{K\times N}$ are chosen such that  $\lambda_j + \mu\lVert \mathbf{W}^T V^j_h\rVert_2^2\geq \max(1,\lVert A^0_h\rVert_2^2) + \epsilon$ for all $\lambda_j \leq \max(1,\lVert A^0_h\rVert_2^2)$, then $\lVert A^k_h\rVert_2,\lVert \phi^k_h\rVert_2\rightarrow 0$.}\end{quote}

It is worth noting that although the bounds above depend on the initial data $\lVert A_h^0\rVert^2_2$ in order to secure the monotone decay of $\lVert A^k_h\rVert_2,\lVert \phi^k_h\rVert_2$, the discrete solution does eventually decay to 0 regardless of $\lVert A_h^0\rVert_2$. Indeed, we have already shown that $A$ can be stabilised with a number of modes that is independent of the initial data. If the bound above is not immediately satisfied, the destabilising term $\lvert A\rvert^2\phi$ in the second equation might be too strong initially, but as $\lvert A\rvert\rightarrow 0$ it is eventually controlled and $\phi$ also begins to decay.

\begin{figure}
\includegraphics[width=.45\linewidth]{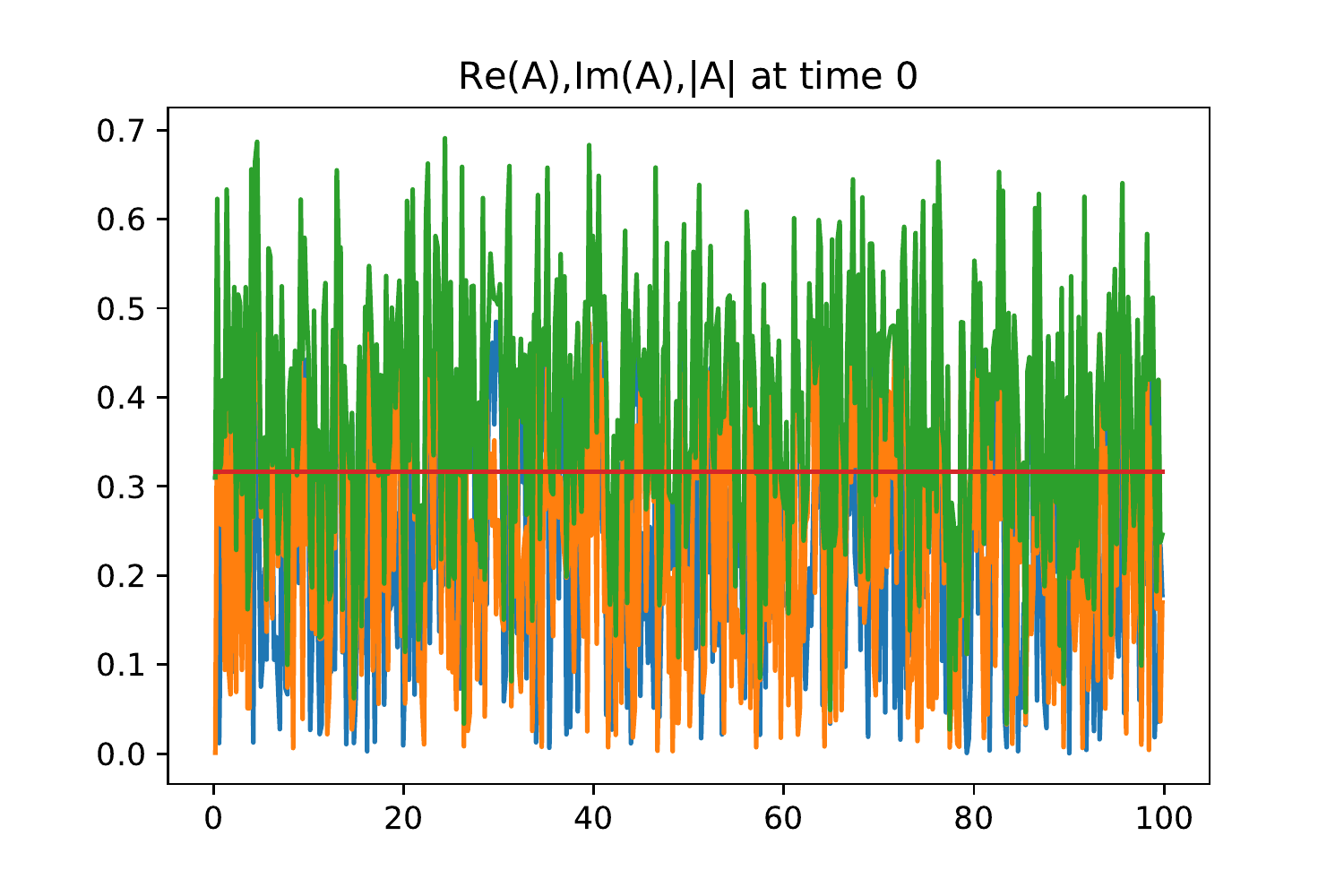}\includegraphics[width=.45\linewidth]{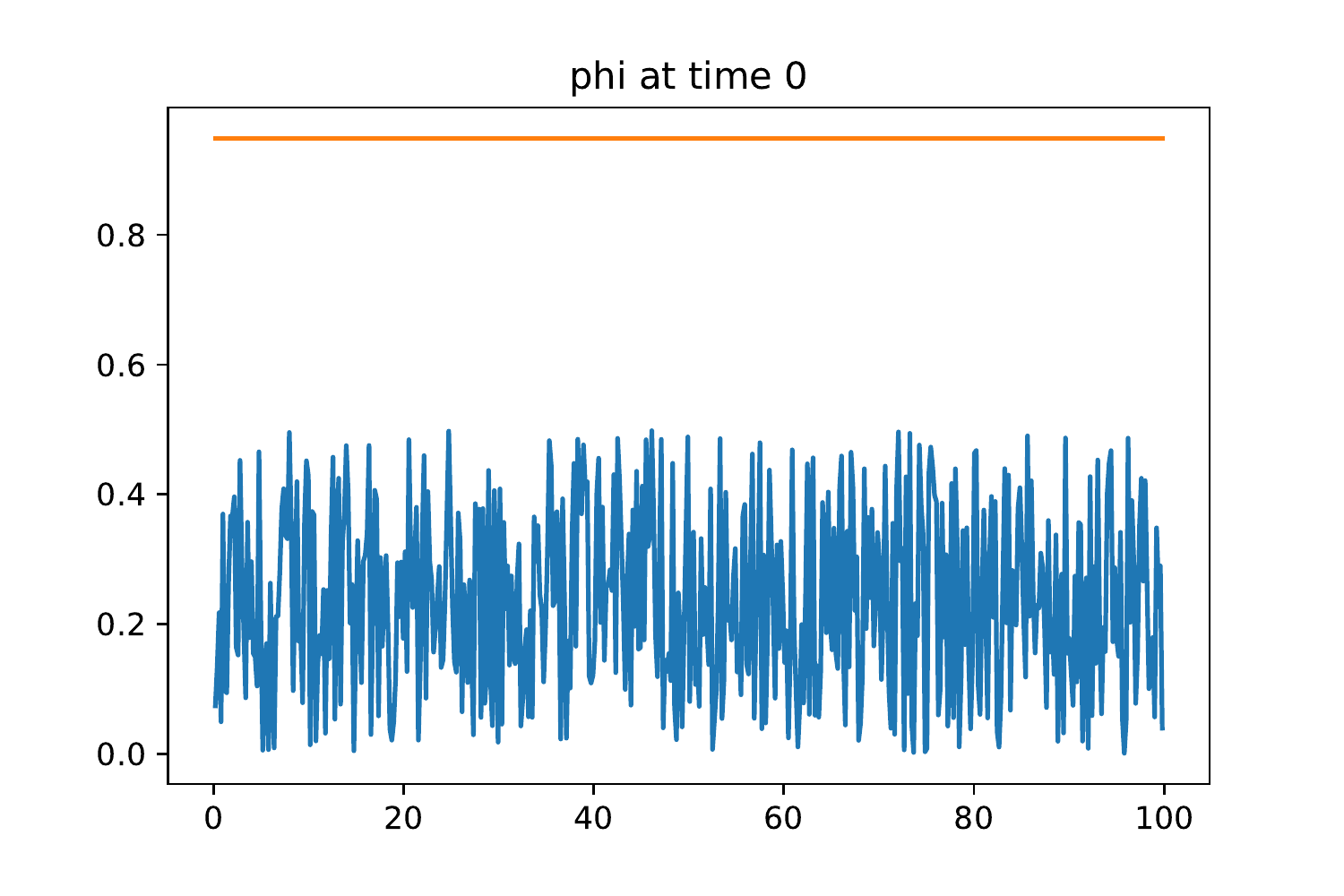}
\includegraphics[width=.45\linewidth]{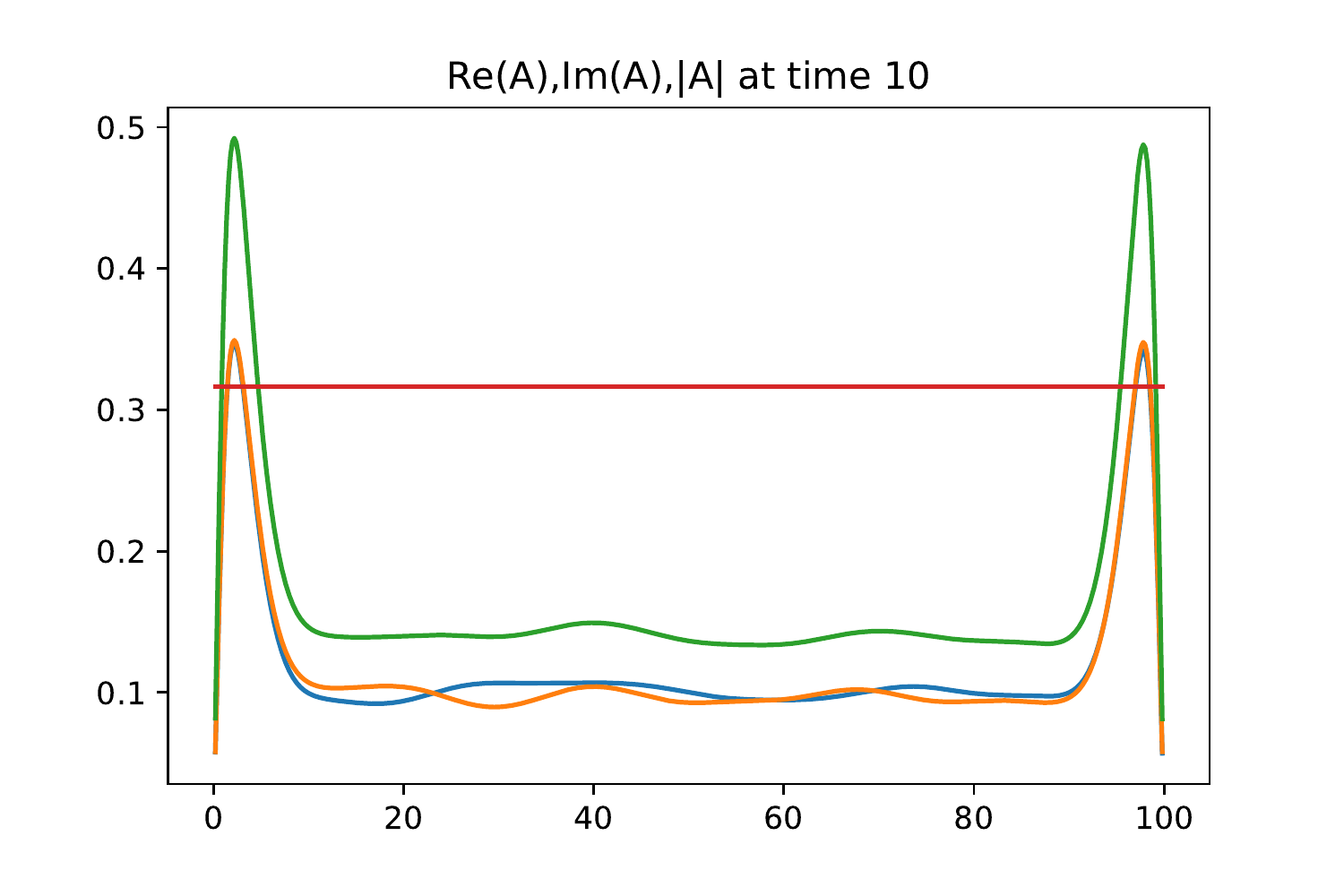}\includegraphics[width=.45\linewidth]{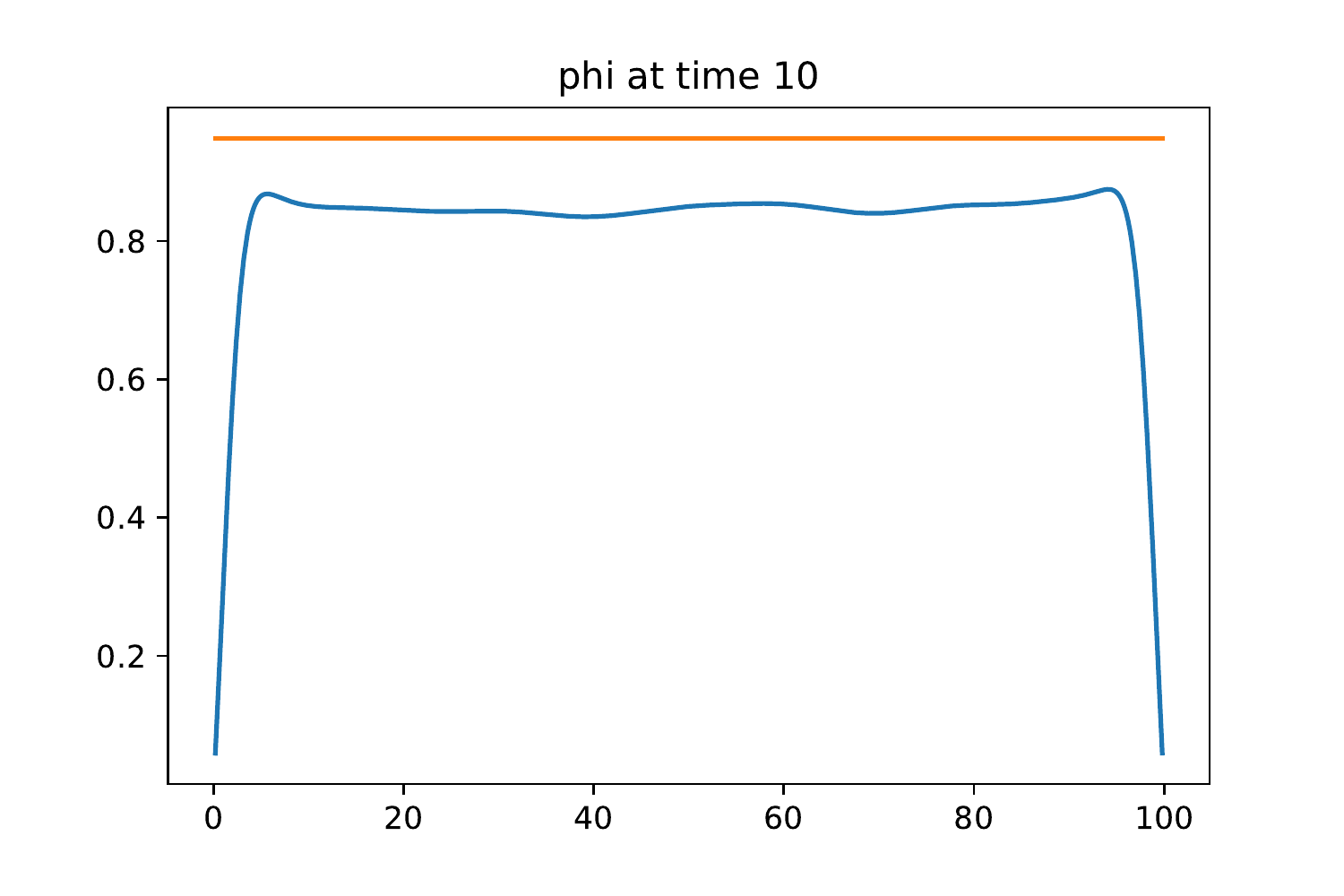}
\includegraphics[width=.45\linewidth]{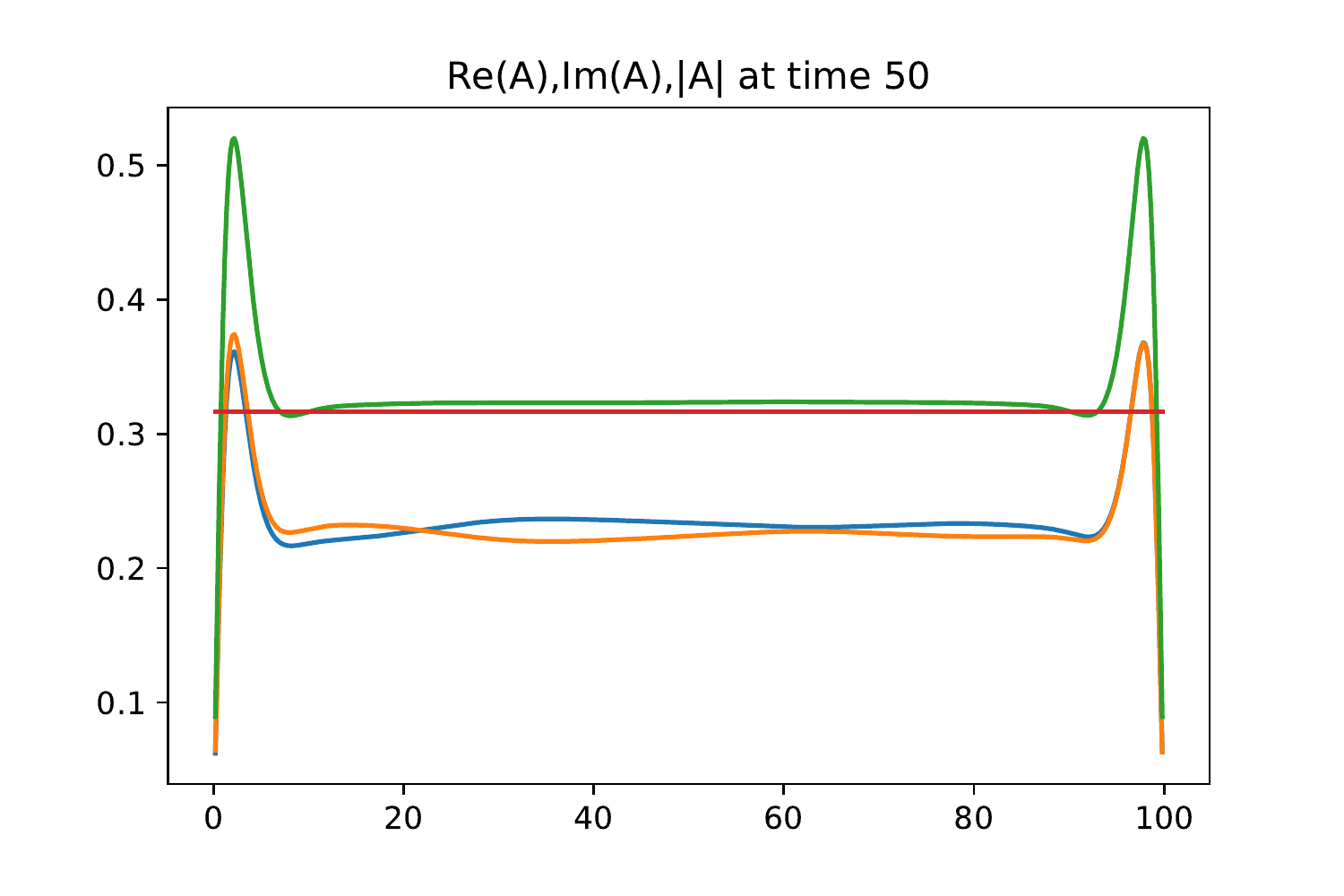}\includegraphics[width=.45\linewidth]{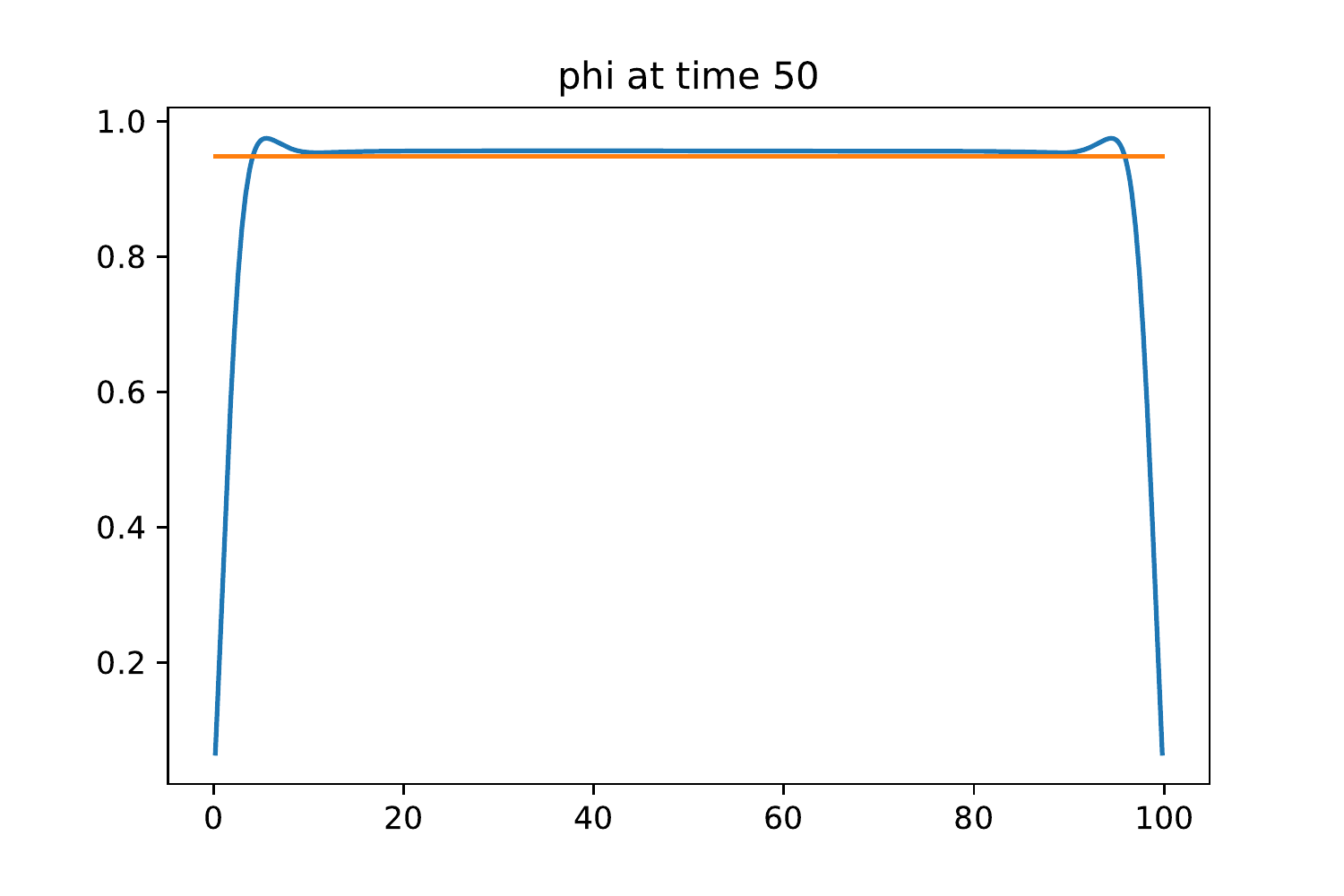}
\includegraphics[width=.45\linewidth]{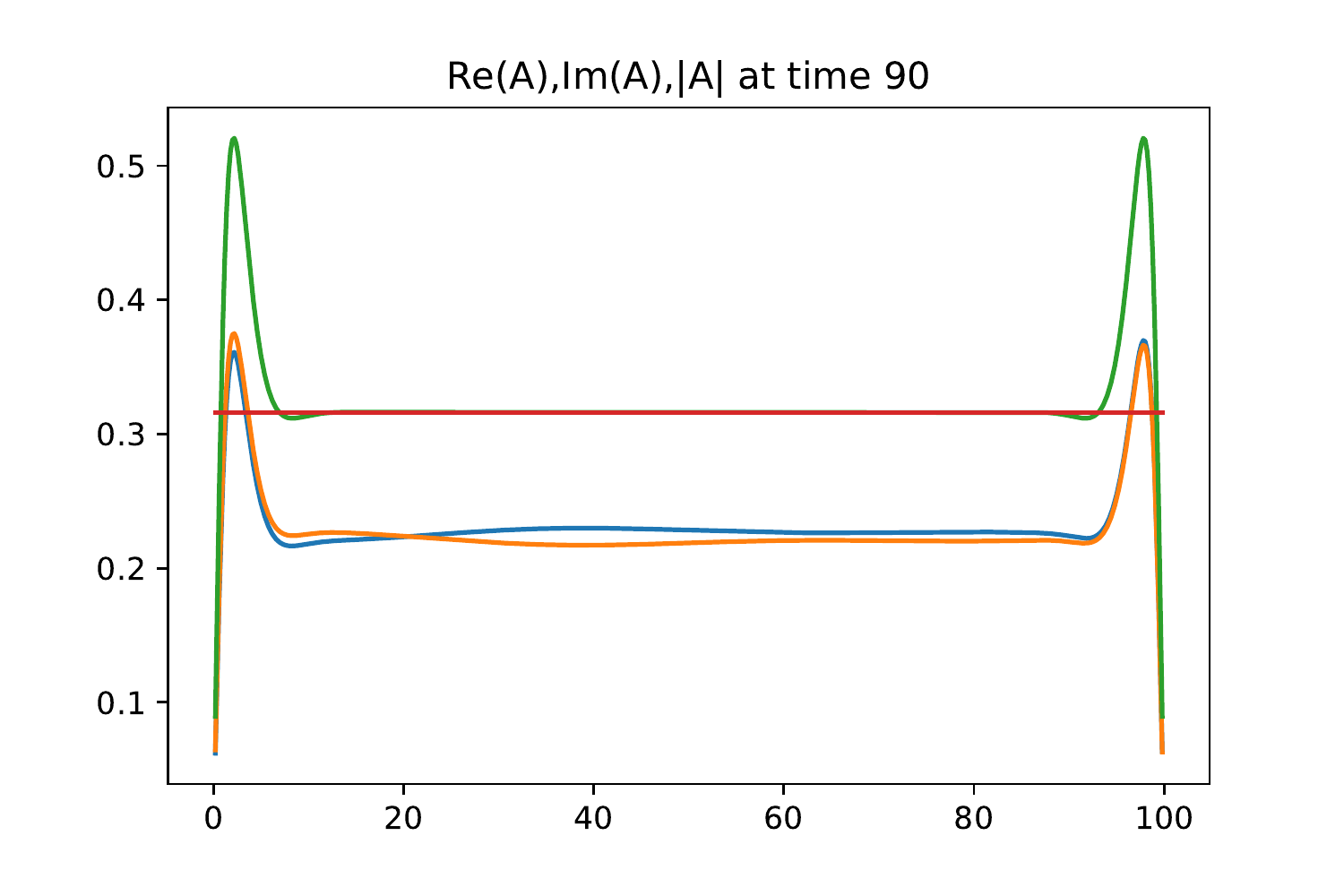}\includegraphics[width=.45\linewidth]{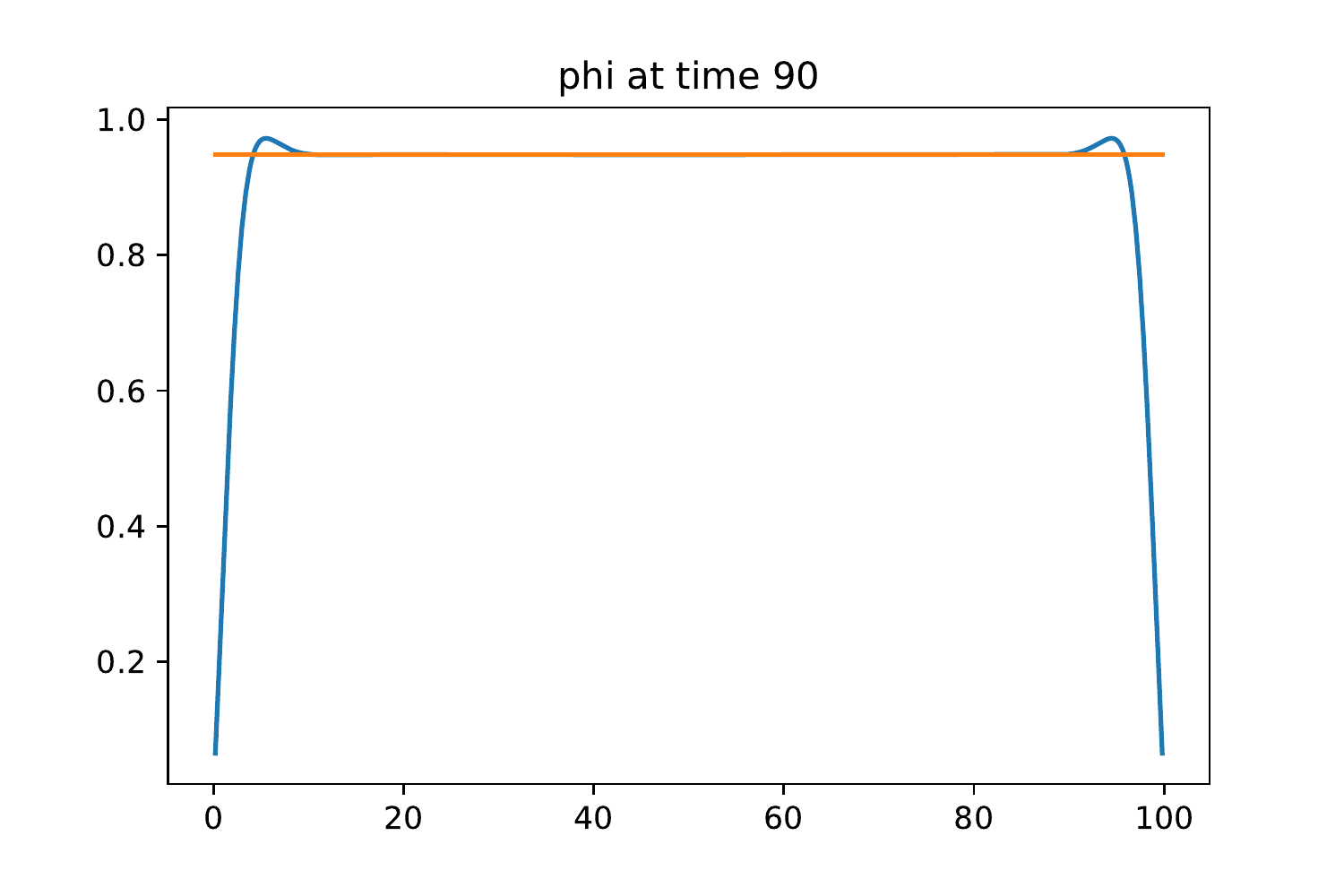}
 \includegraphics[width=.45\linewidth]{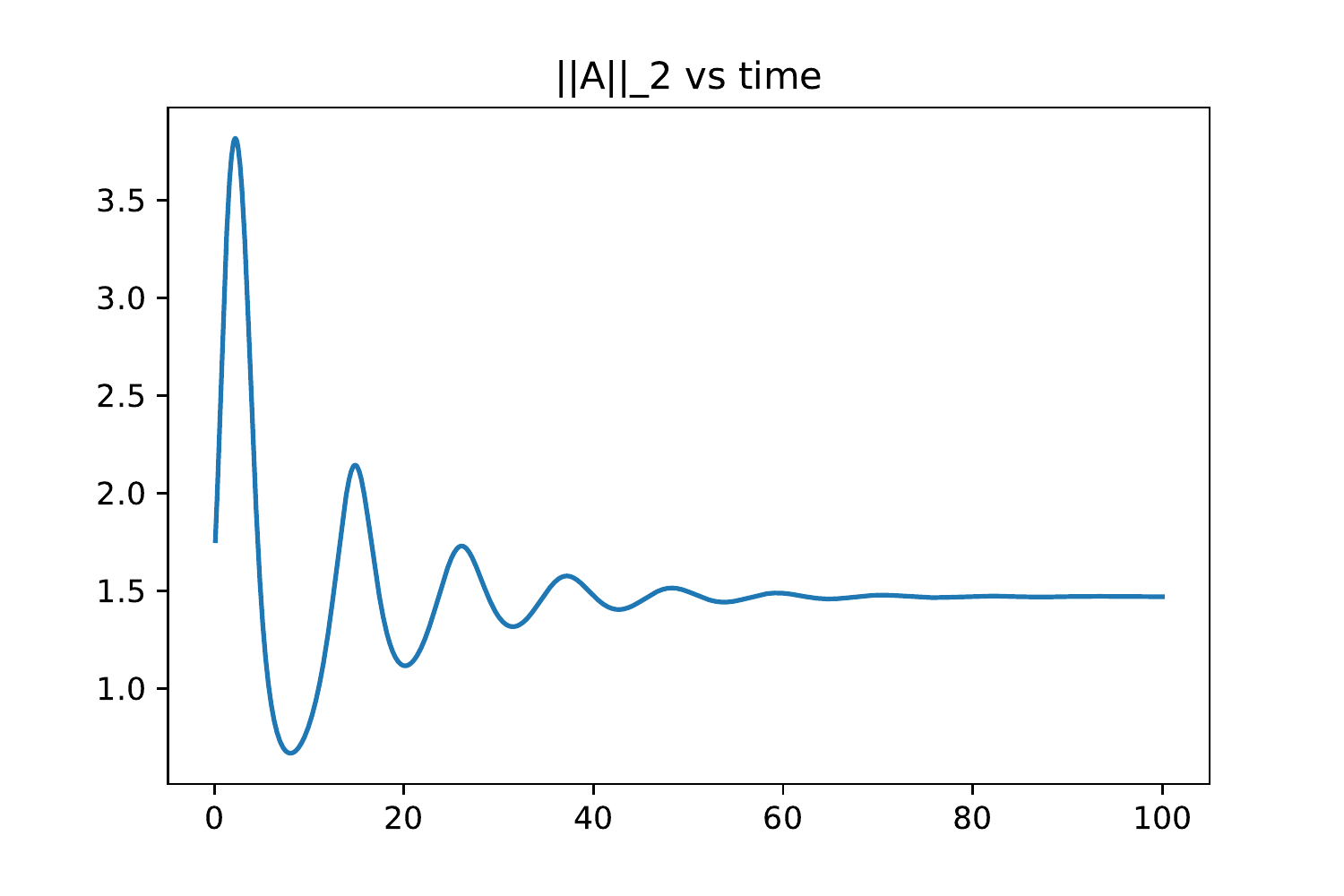}\includegraphics[width=.45\linewidth]{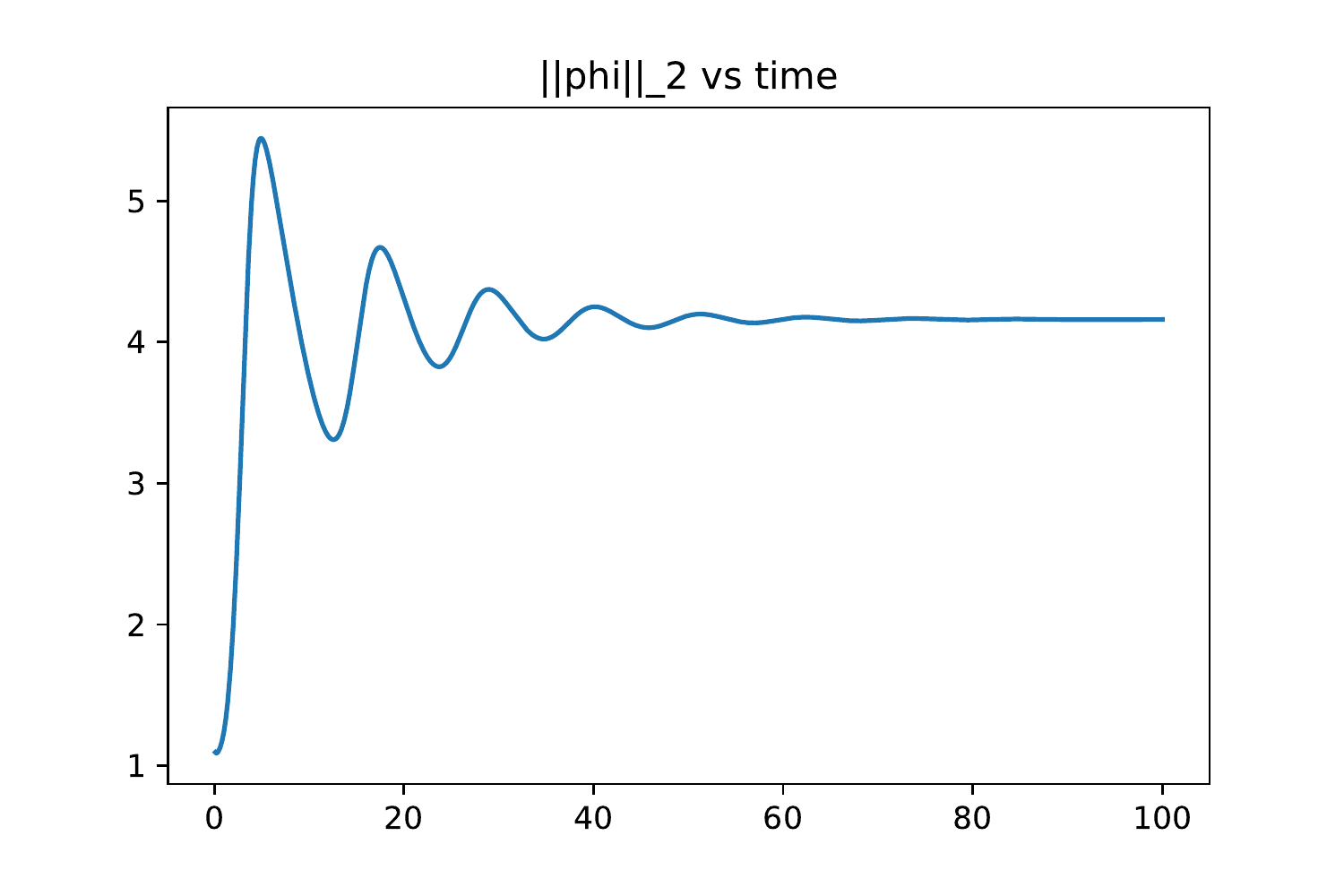}
\caption{Evolution of a highly-oscillatory initial condition towards a non-trivial steady state. The parameters are $\tau = D_1 =1, L=100$ and $h=0.1$.}
\label{fig:nostab}
\end{figure}

\begin{figure}
\includegraphics[width=.45\linewidth]{A0.pdf}\includegraphics[width=.45\linewidth]{phi0.pdf}
\includegraphics[width=.45\linewidth]{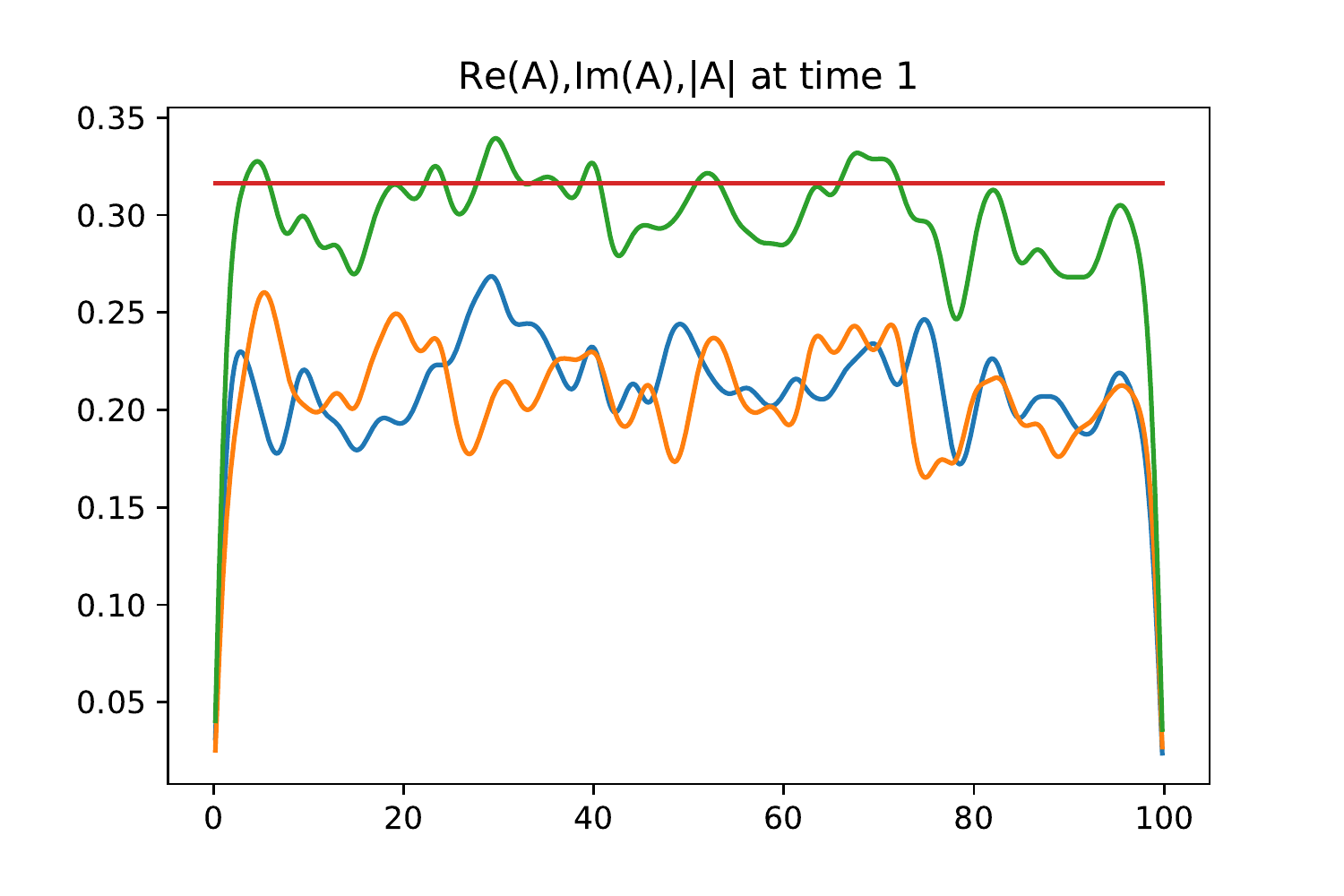}\includegraphics[width=.45\linewidth]{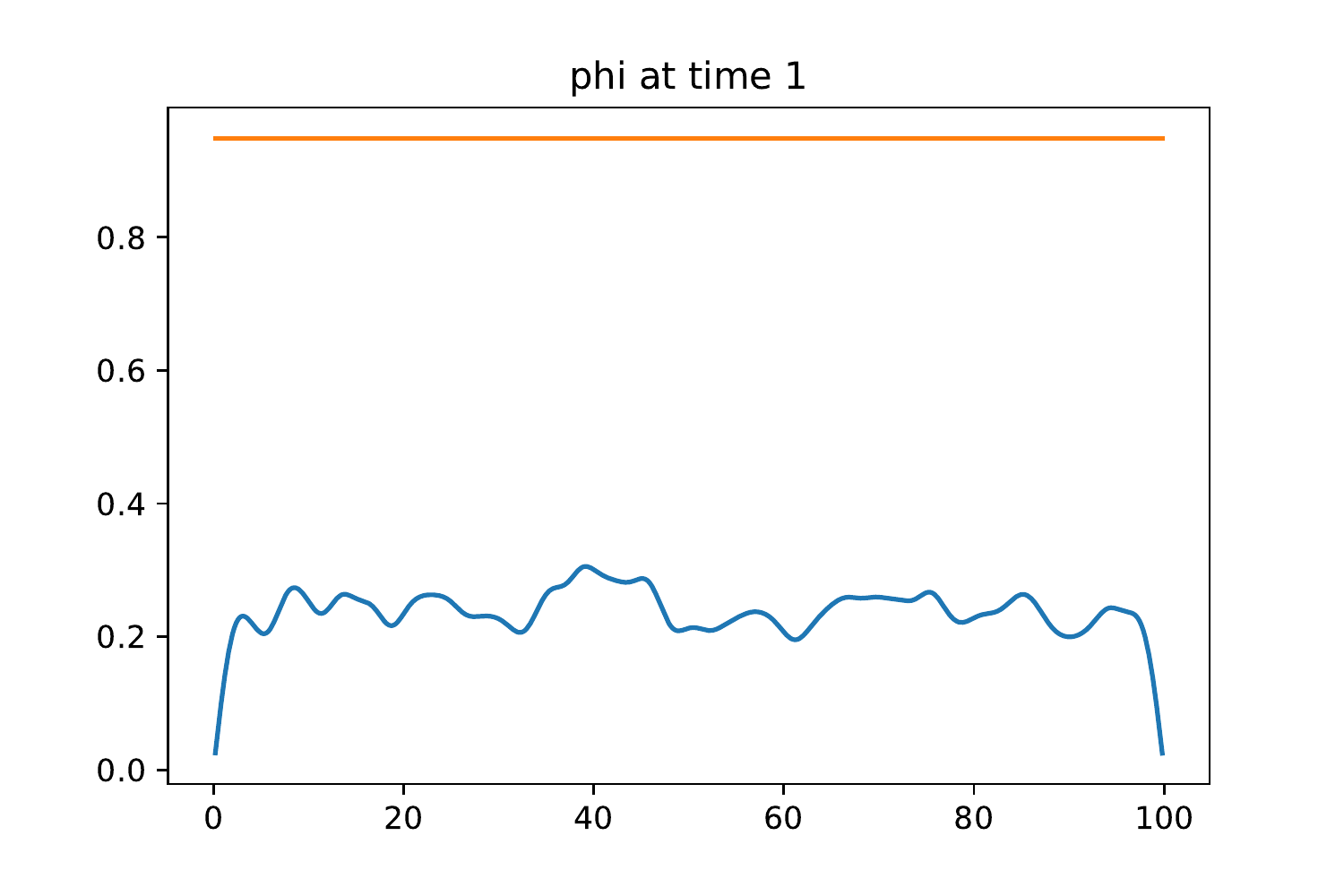}
\includegraphics[width=.45\linewidth]{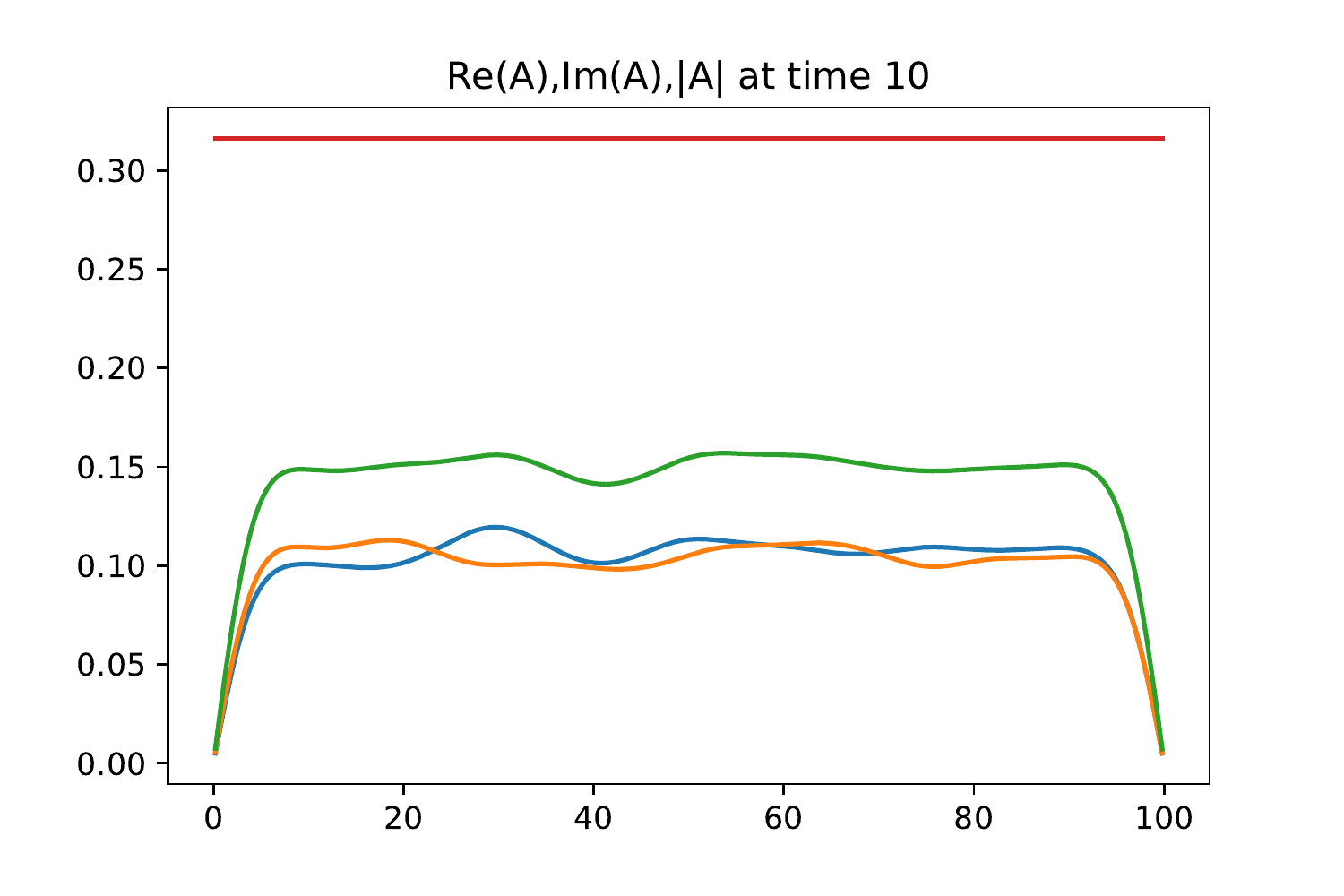}\includegraphics[width=.45\linewidth]{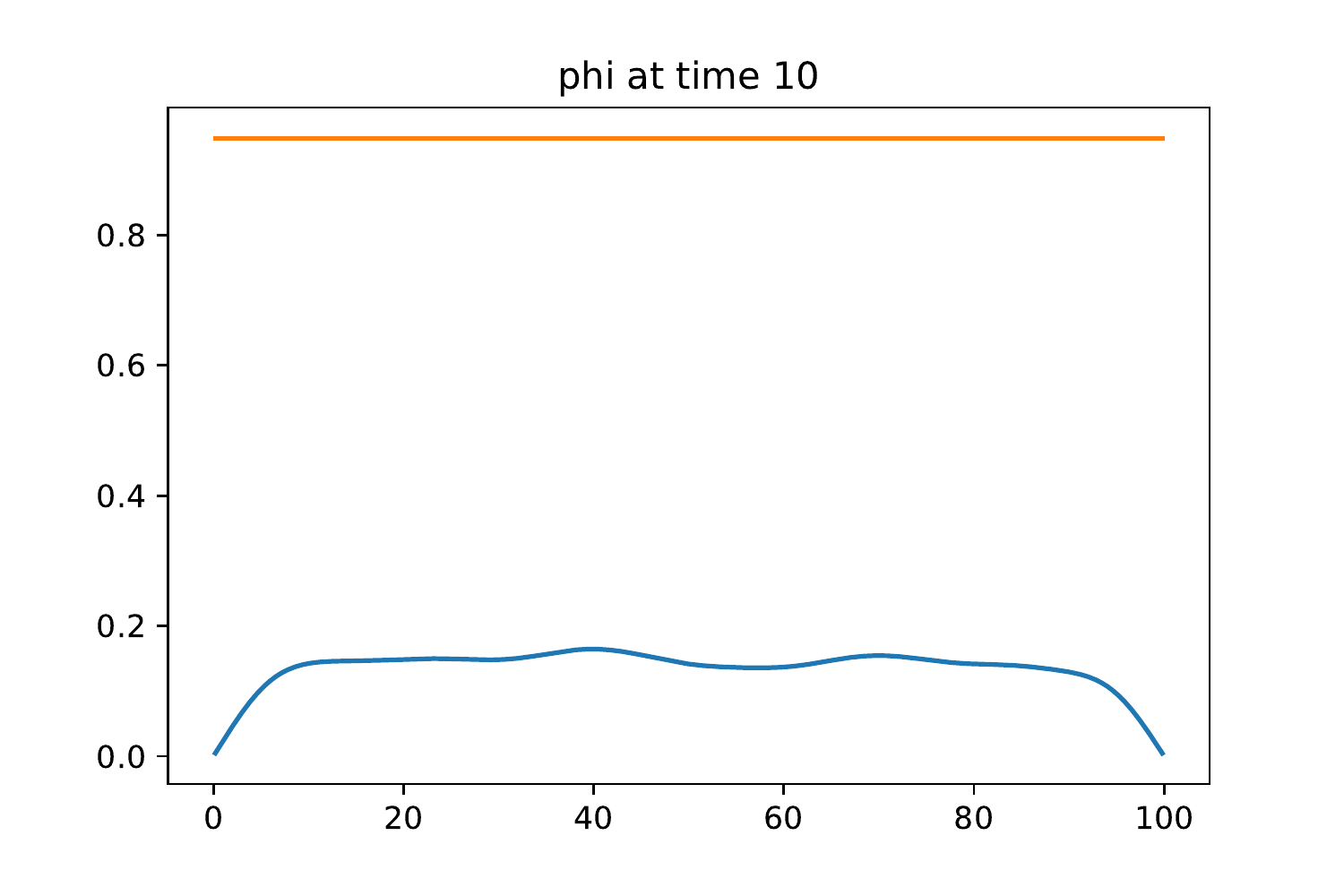}
\includegraphics[width=.45\linewidth]{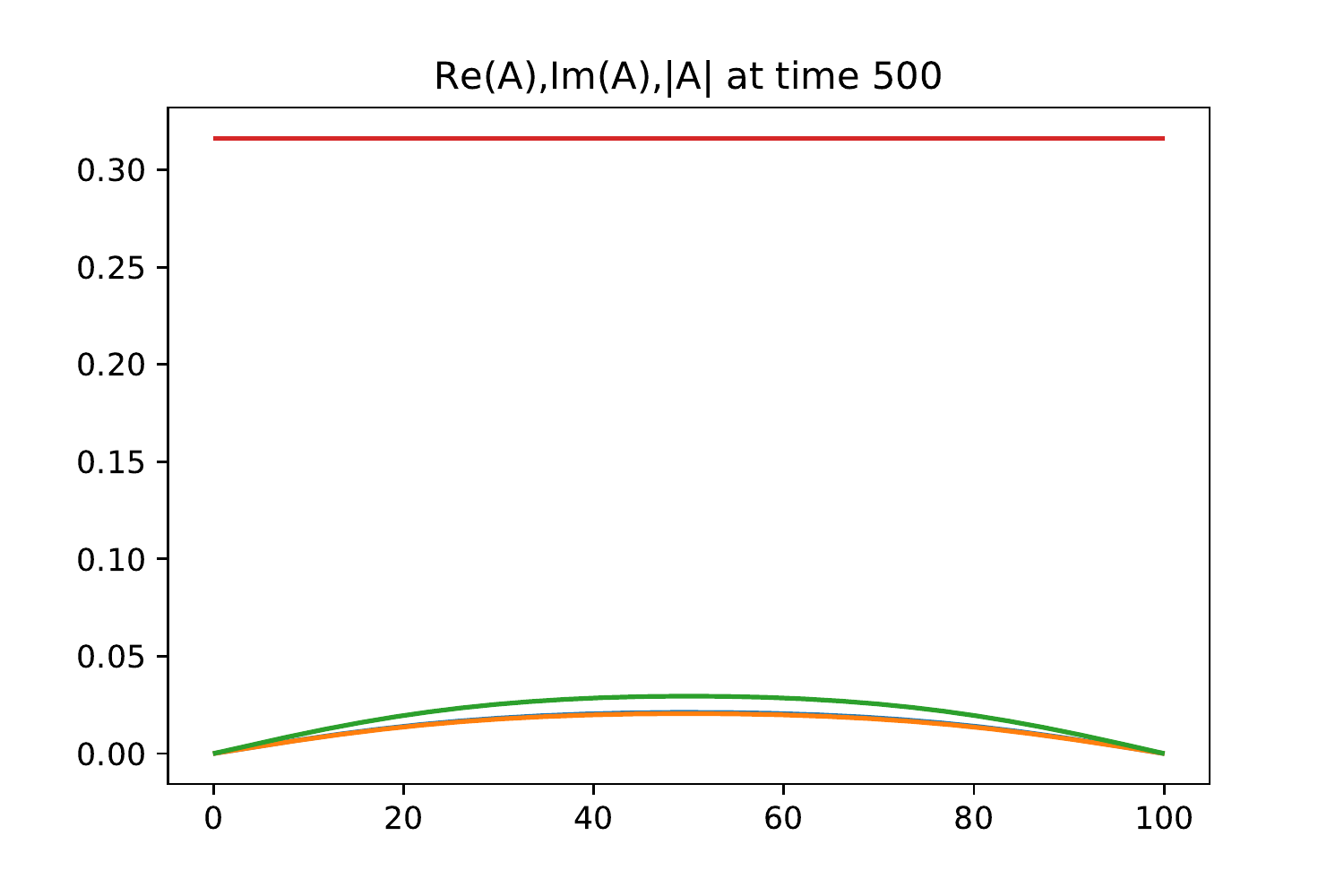}\includegraphics[width=.45\linewidth]{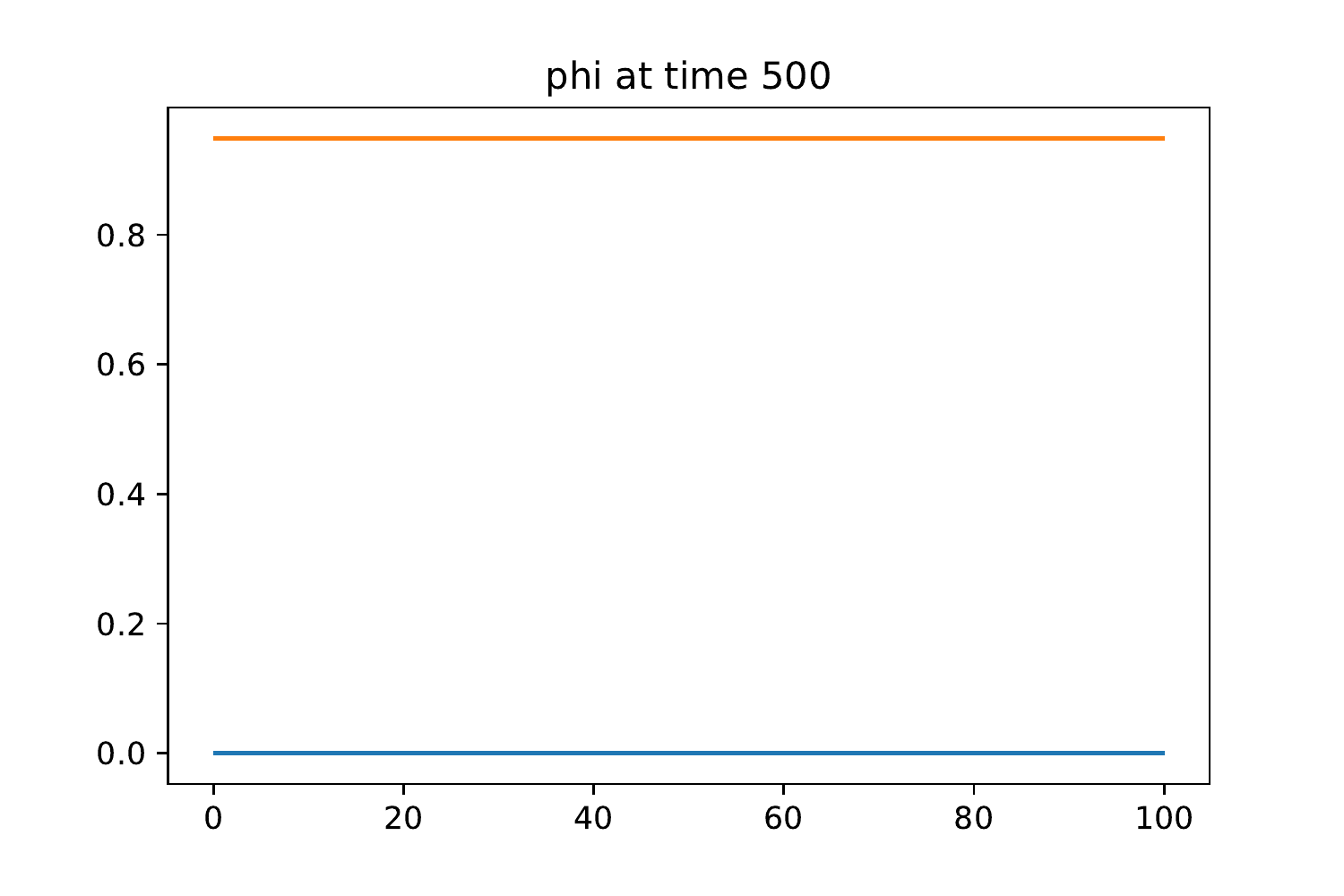}
 \includegraphics[width=.45\linewidth]{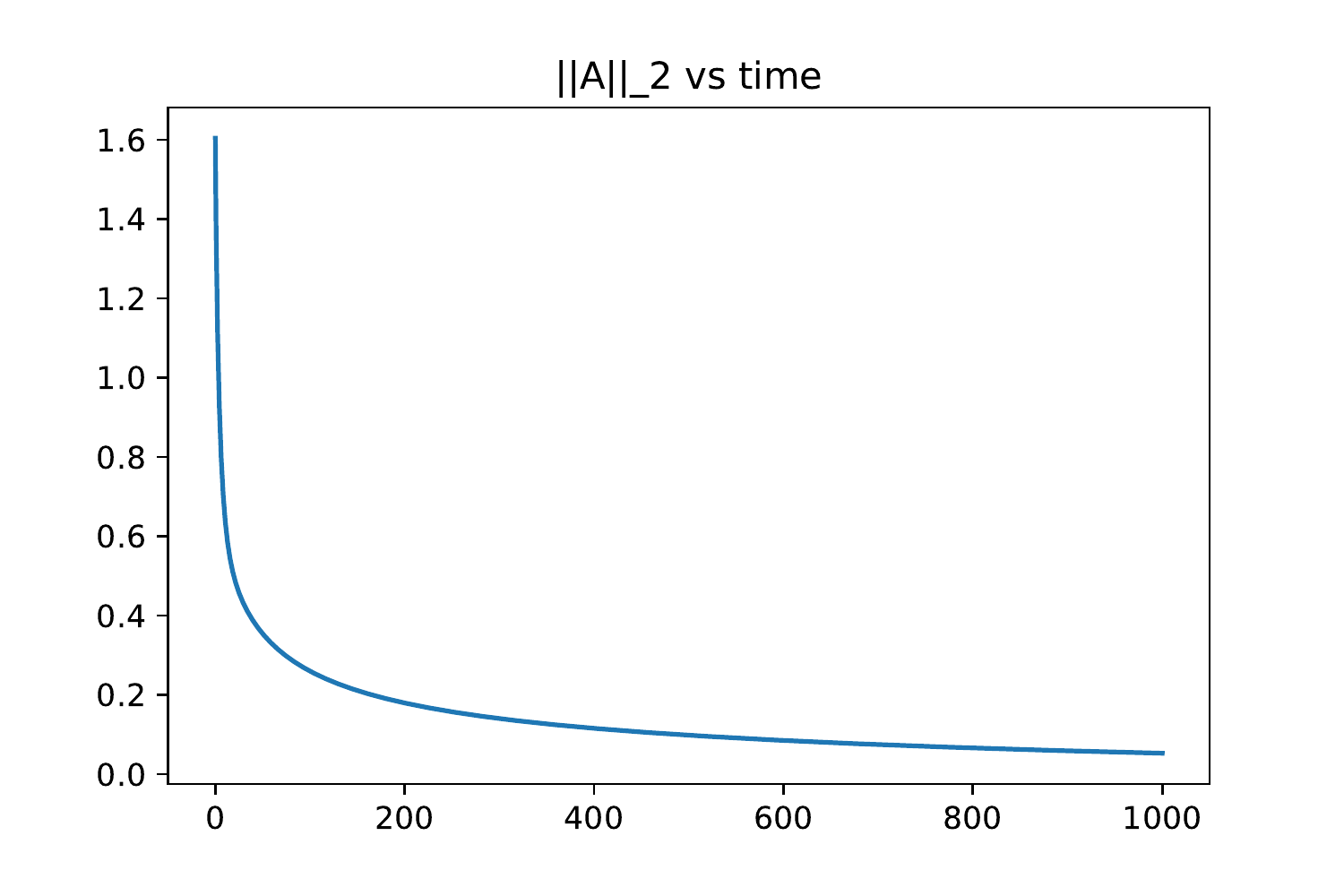}\includegraphics[width=.45\linewidth]{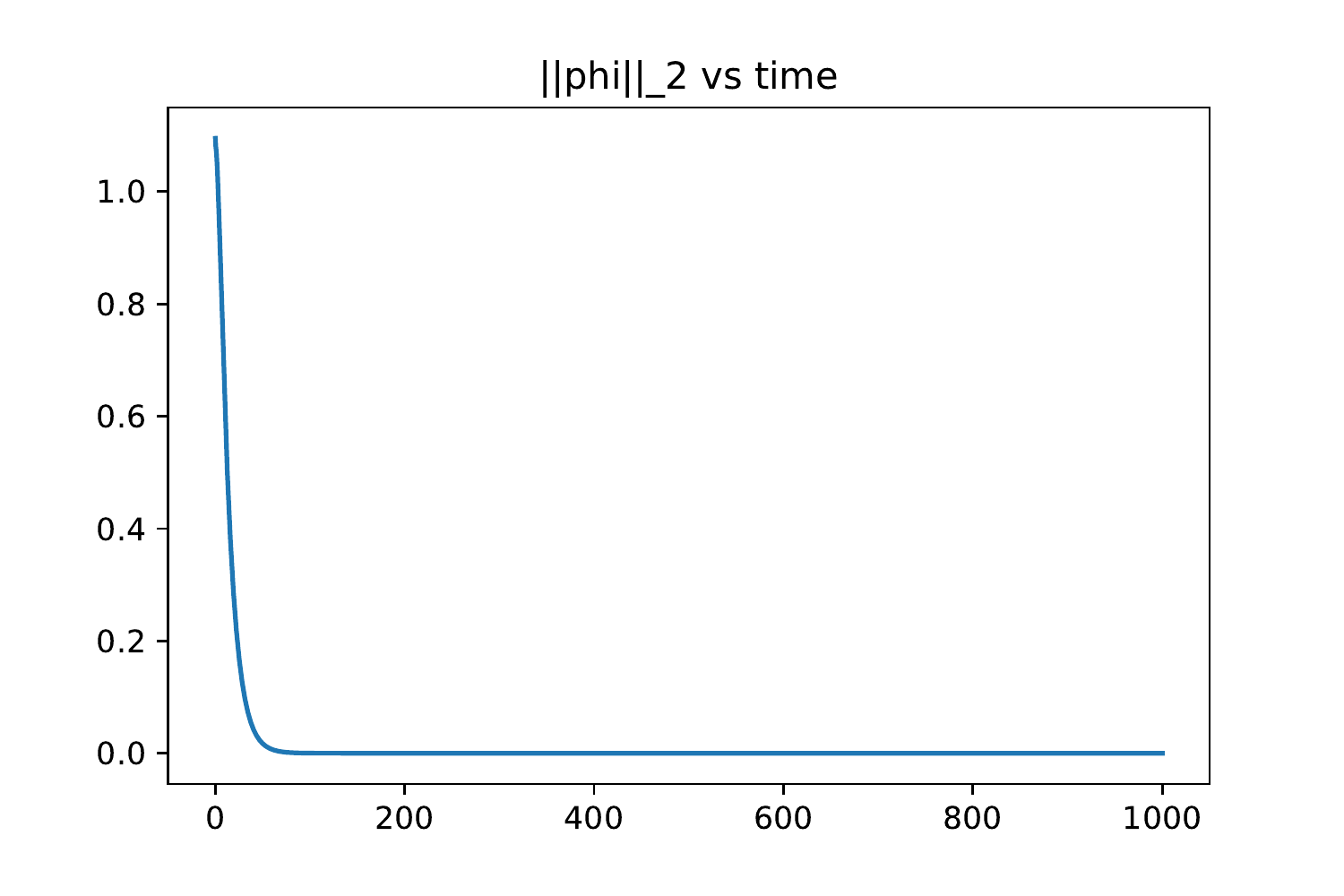}
\caption{Full stabilisation of the solution of Fig.~\ref{fig:nostab} with the $K=147$ Fourier modes prescribed by the stabilization scheme in the text. The solution decays towards the zero state.}
\label{fig:histab}
\end{figure}

Using discrete Fourier modes as $W^j_h$, i.e.~$W^j_i \sim \cos(\frac{j\pi}{L} i \delta x)$ and scaled so that $\lVert V^j_h\rVert_2 =1$, which correspond to the eigenvalues $\lambda_j = (\frac{j\pi}{L})^2$, we can make the choice, albeit somewhat over-cautious, $\mu = \max(1,\lVert A^0_h\rVert_2^2)$ and keep the first $K=\lceil \sqrt{\mu}\frac{L}{\pi}\rceil+1$ modes. In Fig.~\ref{fig:nostab}-\ref{fig:histab}, we present the evolution of a highly oscillatory initial condition, so that all modes are present with high probability, in the presence and absence resp.~of stabilisation. Under no stabilisation, the solution in Fig.~\ref{fig:nostab} tends towards the steady state $\lvert A\rvert = \sqrt{h}$ and $\phi = \sqrt{1-h}$ predicted by the dynamic analysis (modulo boundary layers at the two ends of the domain due to the Dirichlet boundary conditions). Applying full stabilisation with the number of modes prescribed above, we observe in Fig.~\ref{fig:histab} monotone decay of the solution towards the zero state as predicted.

\begin{remark}
Similar to how it is shown in the section 3 we can show that the semigroup $S(t): V^{0}\rightarrow V^{0}$ generated by the problem possesses an exponential attractor.
\end{remark}


\clearpage
    
\bibliography{chevron_pat_stabil.bib}
\bibliographystyle{amsplain}

\end{document}